\documentclass[11pt]{article}

\usepackage{abstract, adjustbox, algorithm, algpseudocode, amsmath, amsthm, amssymb, appendix, bm, booktabs, caption, etoolbox, lettrine, listofitems, listings, lmodern, fancyhdr, graphicx, multirow, neuralnetwork, subfigure, tabularray, tabularx, titlesec, titling, wrapfig, xcolor, xurl}
\usepackage[subfigure]{tocloft}
\usepackage[colorlinks=true, allcolors=blue]{hyperref}
\usepackage[english]{babel}
\usepackage{ebgaramond}
\usepackage[scaled]{helvet}
\usepackage[T1]{fontenc}

\usepackage[margin=0.7in, top=1.2in, bottom=1in]{geometry}

\definecolor{box_color}{HTML}{ECF3F6}
\definecolor{box_color_dark}{HTML}{D0E6F0} 
\definecolor{linkcolor}{HTML}{1B6EA2}

\hypersetup{colorlinks=true,linkcolor=linkcolor,urlcolor=linkcolor,citecolor=linkcolor,allcolors=linkcolor}
\algnewcommand{\Inputs}[1]{%
  \State \textbf{Inputs:}
  \Statex \hspace*{\algorithmicindent}\parbox[t]{.8\linewidth}{\raggedright #1}
}

\pretitle{
    \begin{flushleft}
    \vspace{-15mm} 
    \Huge
    \bfseries 
    \sffamily
} 
\posttitle{
    \end{flushleft}
    \vspace{3mm}
    \hrule
    \vspace{3mm}
} 
\preauthor{\begin{flushleft}\large\sffamily}
\postauthor{\end{flushleft}}

\renewenvironment{abstract}
 {\par\noindent\textbf{\sffamily \abstractname} \ \ignorespaces}
 {\par\medskip}
 
\renewcommand\thesection{\Roman{section}} 
\renewcommand\thesubsection{\roman{subsection}} 
\titleformat{\section}[block]{\Large\sffamily\bfseries}{\thesection.}{1em}{} 
\titleformat{\subsection}[block]{\large\bfseries\sffamily}{\thesubsection.}{1em}{} 

\DeclareCaptionFormat{caption_format}{\normalfont\sffamily\fontsize{8}{10}\selectfont#1#2#3}
\adjustboxset{margin=10pt,bgcolor={HTML}{ECF3F6},rndcorners=6}

\DeclareRobustCommand{\tablefont}{%
        \fontencoding{\encodingdefault}%
        \fontseries{m}
        \fontshape{n}
        \fontfamily{phv}
        \fontsize{9}{12}
        \selectfont}
\DeclareTextFontCommand{\texttable}{\tablefont}

\newcommand{\dropcap}[1]{\lettrine[lines=2,lraise=0.05,findent=0.1em, nindent=0em]{{\sffamily{#1}}}{}} 
\newcommand{\cs}[1]{\sffamily\fontsize{8}{10}\selectfont{#1}} 

\newcommand{\predm}[1]{\bm{\hat{#1}}}

\usepackage{csquotes}
\usepackage[backend=biber,style=nature]{biblatex}
\AtEveryBibitem{%
  \ifboolexpr{
    test {\ifentrytype{article}} 
    or
    test {\ifentrytype{book}}
  }{%
    \clearfield{url}%
    \clearfield{urldate}%
    \clearfield{urlyear}
    \clearfield{month}
    \clearfield{issn}
  }{}%
}
\usepackage{ragged2e}
\AtBeginBibliography{\RaggedRight}

\newif\ifhighlight

\highlightfalse

\definecolor{highlightcolor}{HTML}{75B09C}

\newcommand{\highlight}[1]{{\ifhighlight\texorpdfstring{\color{highlightcolor}#1}{#1}\else#1\fi}}

\title{Modelling Global Trade \\ with Optimal Transport} 

\date{}

\author{
	\bfseries{Thomas Gaskin}$^{\text{1,2}\star}$,
     \bfseries{Guven Demirel$^\text{3}$},
     \bfseries{Marie-Therese Wolfram}$^\text{4}$, \bfseries{Andrew Duncan}$^\text{2}$
}

\begin{document}

\maketitle

\vspace{-20mm}
{\small\flushleft
\textsuperscript{\textbf{1}}~Department of Applied Mathematics and Theoretical Physics, University of Cambridge, Cambridge CB3 0WA, United Kingdom; \textsuperscript{\textbf{2}}~Department of Mathematics, Imperial College London, London SW7 2AZ, United Kingdom; \textsuperscript{\textbf{3}}~School of Business and Management, Queen Mary University of London, London E1 4NS, United Kingdom; \textsuperscript{\textbf{4}}~Mathematics Institute, University of Warwick, Coventry CV4 7AL, United Kingdom

\medskip 
Correspondence: $^\star$trg34@cam.ac.uk

}

\vspace{6mm}

\hrule
\vspace{6mm}
\begin{abstract} 
\noindent Global trade is shaped by a complex mix of factors beyond supply and demand, including tangible variables like transport costs and tariffs, as well as less quantifiable influences such as political and economic relations. Traditionally, economists model trade using gravity models, which rely on explicit covariates \highlight{that might} struggle to capture these subtler drivers of trade. In this work, we employ optimal transport and a deep neural network to learn a time-dependent cost function from data, without imposing a specific functional form. This approach consistently outperforms traditional gravity models in accuracy and \highlight{has similar performance to three-way gravity models}, while providing natural uncertainty quantification. Applying our framework to global food and agricultural trade, we show that the Global South suffered disproportionately from the war in Ukraine's impact on wheat markets. We also analyse the effects of free-trade agreements and trade disputes with China, as well as Brexit's impact on British trade with Europe, uncovering hidden patterns that trade volumes alone cannot reveal.
\noindent 

\bigskip

\end{abstract}
\vspace{6mm}
\hrule
\vspace{1mm}
\setcounter{tocdepth}{1}
\tableofcontents
\vspace{4mm}

\newpage

\twocolumn
\section{Introduction}
\vspace{2cm}

\dropcap{I}nternational trade serves as the backbone of the world economy, distributing goods and connecting markets through global logistics networks. Its dynamics are driven by numerous factors beyond mere supply and demand, such as tariffs, non-tariff policy barriers, political and economic tensions, and disruptions caused by accidents, conflicts, and civil wars. Among all traded commodities, agricultural and food products hold particular interest for policymakers and the general public due to their significant volume, high trade value, and critical role in food security and resilience \cite{Friel_et_al2020, Wood_et_al2023}. Consumer food prices are a product of all the complexly interwoven factors governing trade. However, they do not always directly reflect the ease of doing business between any two countries. For instance, in May 2020, China imposed an 80\% tariff on Australian barley, leading to a major restructuring of global supply chains (see fig.~\ref{fig:Barley_trade}): Chinese demand was suddenly met from France, Canada, and Argentina, while Australia started exporting surplus barley e.g. to Saudi Arabia. Despite these shifts, for the next five months the global barley price barely budged \cite{GIWA_2020, FRED_Barley}.

Modelling global trade has garnered significant attention in the economics literature, with \emph{gravity models} being the most widely used approach \cite{Eaton_Kortum_2002, Anderson_Wincoop_2003, Arkolakis_2012, Allen_2014, Yotov2022}. These models, named for their direct analogy to Newton's law of gravity, assume that the total trade $T_{ij}$ of a given commodity between two countries $i$ and $j$ is proportional to the total output $O_i$ of the source country and the total expenditure $E_j$ of the destination country, as well as being inversely related to a `distance' between them:
\begin{equation}
    T_{ij}(t) \sim \dfrac{O_i(t) E_j(t)}{C_{ij}(t)}.
    \label{eq:Gravity}
\end{equation}
This distance $C_{ij}$ comprises all factors that contribute to the ease of selling goods produced in one country to another, including transportation costs, information costs, and tariff and non-tariff barriers to trade. Traditional gravity models use a set of covariates to estimate $C_{ij}$ as
\begin{align}
    \log C_{ij}(t) = & \sum_k \alpha_k \pi_{i, k}(t) + \sum_l \beta_l \chi_{j, l}(t) \nonumber \\ & + \sum_m \gamma_m \rho_{ij, m}(t),
    \label{eq:Gravity_covariates}
\end{align}
where $\pi_i$ and $\chi_j$ are exogenous exporter and importer-side regressors \cite{Dekle_2008}, \highlight{which can be specified as fixed effects}, $\rho_{ij}$ are bilateral covariates, and $\alpha$, $\beta$, $\gamma$ are the coefficient vectors. Commonly used covariates include geographic proximity, the existence of trade agreements, colonial ties, tariffs, non-tariff barriers, or shared languages \cite{Yotov_2017}. The \emph{structural gravity model} corrects eq.~\eqref{eq:Gravity} with import and export multilateral resistance terms, which account for the relative nature of bilateral trade shares. This adjustment has been shown to align with various microeconomic models \cite{Arkolakis_2012}. Gravity models have been widely used to study agrifood trade. For instance, \cite{Olper_Raimondi2009} estimate residual trade costs based on a micro-founded gravity equation, finding ad-valorem costs to be 60\% higher in the Global South compared to the North. Studies have also investigated the impact of global and regional trade agreements \cite{Sarker_Jayasinghe2007,Mujahid_Kalkuhl2016} and the effect of eliminating tariffs \cite{Raimondi_Olper2011,Philippidis_etal2013}. 

The gravity-based approach is attractive to researchers due to its interpretability, mathematical simplicity, and consistency with various microeconomic theories \cite{Yotov2022}. However, it is not without its limitations. For one, multilateral trade resistance terms, central to the structural gravity model, are unobservable and must be estimated, often using fixed effects \cite{Yotov_2017}. Elasticity and other key parameters are often unavailable at a granular level, requiring aggregation that can introduce bias \cite{Breinlich2022}. The model's cost function also depends heavily on the choice of covariates and functional form, making specification crucial for interpreting results. In addition, unobservables---such as the subtle effects of changing political relations, public preferences, or aversions toward products from specific countries---are absorbed in the error term. Finally, while trade costs are generally asymmetric (\(C_{ij} \neq C_{ji}\)), commonly used covariates are not, making it difficult for a model to capture the inherent imbalances in trade relationships. See \cite{Yotov_2017, Yotov2022, Capoani2023} for a deeper discussion of challenges and best practices.

In this work, we present a more general approach that dispenses with the use of covariates and a functional form, instead inferring the cost directly from data. Our method is based on the \emph{optimal transport} (OT) framework \cite{Villani_2021}, which generalises gravity-based models. In OT, trade flows are assumed to match supply and demand to minimise an overall cost. Mathematically, this is expressed as follows: let $\bm{C} \in \mathbb{R}^{m \times n}$ be a matrix quantifying the `cost' (in a general sense) of moving goods from country $i$ to $j$. Given the supply vector $\bm{\mu} \in \mathbb{R}^m$ and the demand vector $\bm{\nu} \in \mathbb{R}^n$, the optimal transport problem consists in finding a \emph{transport plan}, i.e. a matrix $\bm{T} \in \mathbb{R}^{m \times n}_+$ with entries $T_{ij}$ modelling the total volume (or value) of transport from country $i$ to $j$, such that the total cost 
\begin{equation}
	c(\bm{T}) = \sum_{i,j} T_{ij}C_{ij}
	\label{eq:OT_cost}
\end{equation}
is minimised. In addition, the \emph{marginal constraints} 
\begin{equation}
	\sum_i T_{ij} = \bm{\mu}, \ \sum_i T_{ij} = \bm{\nu}
	\label{eq:OT_constraints}
\end{equation}
must be satisfied, ensuring that demand and supply are met. It is advantageous to add a \emph{regularisation term} term to the cost, as it ensures existence of a unique solution and significantly improves computational efficiency; the total cost then becomes
\begin{equation}
\label{eq:cost_eps}
    c_{\varepsilon}(\bm{T}) = c(\bm{T}) \highlight{-} \varepsilon \mathcal{H}(\bm{T}),
\end{equation}
where \highlight{$\mathcal{H}(\bm{T}) = -\sum_{ij} T_{ij} (\log T_{ij} - 1)$ denotes the negative entropy }of $\bm{T}$ and $\varepsilon > 0$ is a regularisation parameter. It can be shown that the solution will then be of the form
\begin{equation}
	\bm{T} = \bm{\Pi} e^{-\bm{C} / \varepsilon} \bm{\Omega},
	\label{eq:OT}
\end{equation}
where $\bm{\Pi}$ and $\bm{\Omega}$ are diagonal scaling matrices which ensure that the marginal constraints hold (see Methods). As described in \cite{Wilson_1967}, gravity models can be reformulated as solutions of a regularised OT problem with an appropriate choice of parameters. While OT-based models might appear to suggest a centralised control of flows, its dual formulation admits an alternative, decentralised interpretation of importers seeking to minimise the cost of purchasing commodities and exporters seeking to maximise their profit (see Methods). The solution at equilibrium coincides with the solution of the OT problem \cite{Galichon_2016}, which in its classic form \eqref{eq:OT_cost}--\eqref{eq:OT_constraints} is well understood. This is less true for the corresponding \emph{inverse problem} we are interested in, despite its mathematical and practical importance: given a (possibly noisy) observation of $\bm{T}$, $\bm{\mu}$, and $\bm{\nu}$, this problem consists in inferring the underlying cost $\bm{C}$. \highlight{Maximum likelihood estimation of costs relates to the inverse optimal transport problem; however, the specific parametrization of the cost matrix in gravity models demands careful estimator design. Moreover, zeros and heteroscedasticity in observed trade flows cause misspecification in gravity model estimation, affecting estimator performance (see \cite{silva2006log}).}

\begin{figure*}[htp!]
\begin{adjustbox}{minipage=\textwidth-20pt}
\begin{minipage}{\textwidth}
	\cs{\textbf{A} Ukrainian wheat exports in metric tons, 2021 (left) and 2022}
\end{minipage}
\hfill \vspace{-3mm}

\begin{minipage}{\textwidth}
	\includegraphics[width=\textwidth]{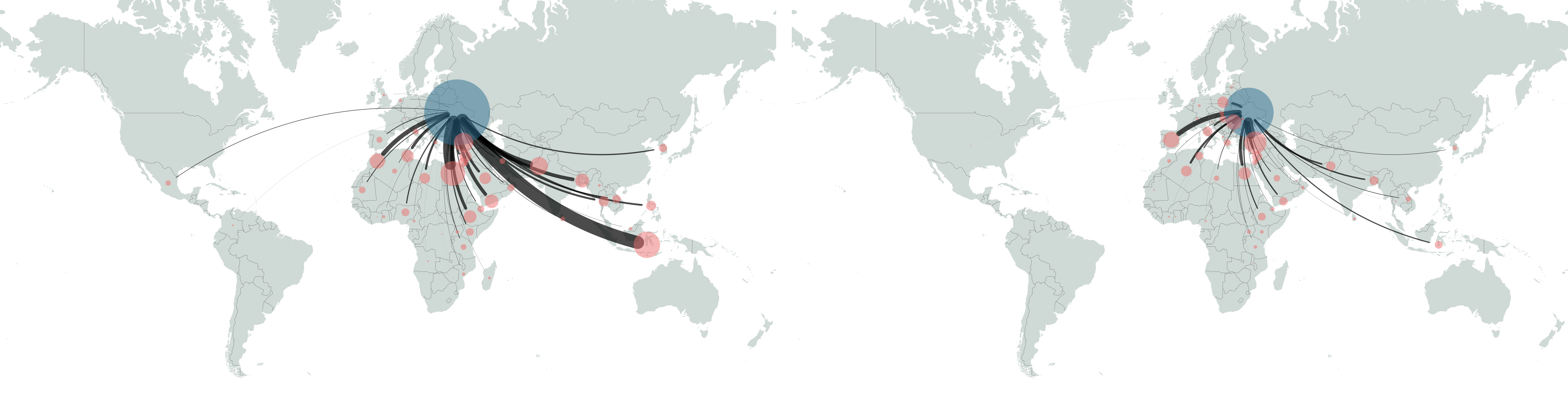}
\end{minipage}
\vspace{1mm}

\begin{minipage}{\textwidth}
	\flushleft \cs{\textbf{B} \highlight{Percent change in trade volume (left) and absolute change in cost}}
\end{minipage}
\hfill \vspace{-3mm}

\begin{minipage}{\textwidth}
	\includegraphics[width=\textwidth]{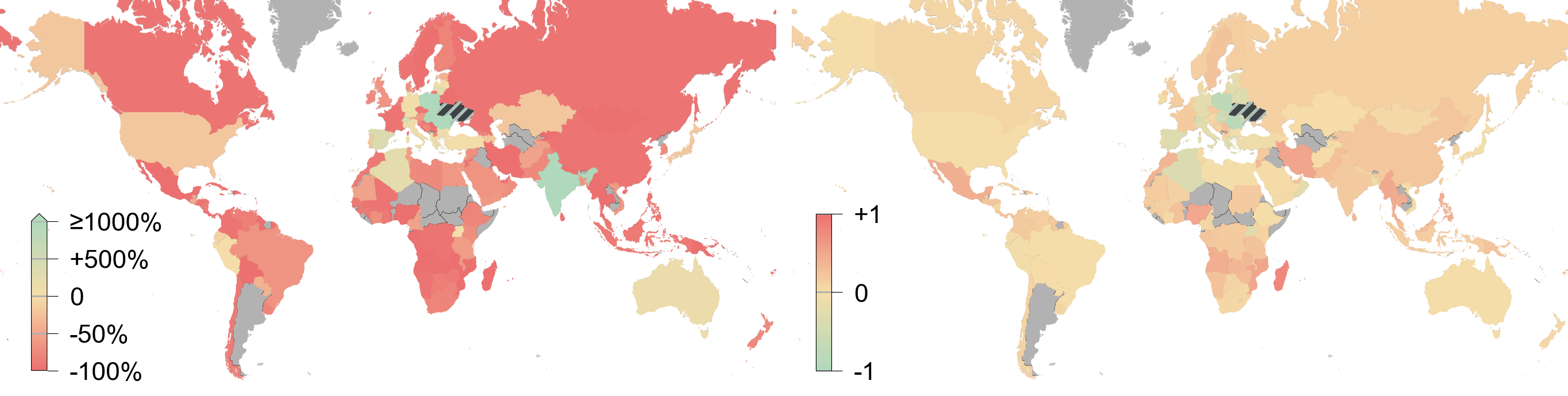}
\end{minipage}
\vspace{1mm}

\begin{minipage}{\textwidth}
	\cs{\textbf{C} \highlight{Change in trade volume and cost, selected countries}}
\end{minipage}
\hfill \vspace{-3mm}

\begin{minipage}{\textwidth}
	\includegraphics[width=\textwidth]{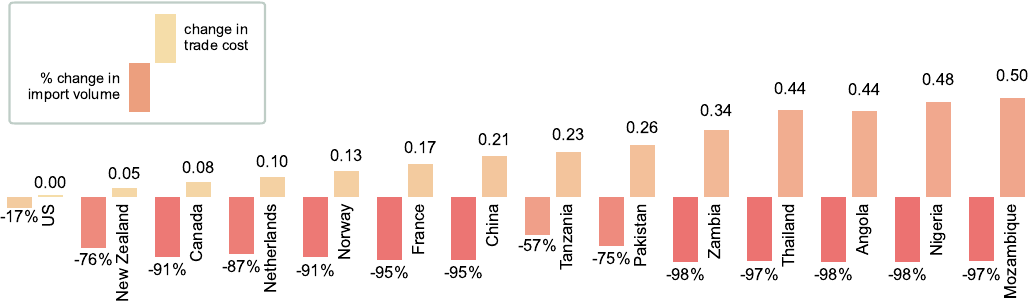}
\end{minipage}

\caption{Ukrainian wheat exports, 2021--2022. \textbf{A} Network of Ukrainian exports, 2021 and 2022. Shown are the largest trading partners, making up 99\% of Ukrainian exports. The blue node represents the total Ukrainian export volume (in metric tons), the red nodes are the import volumes. Edge widths represent the flow volume. \textbf{B} The change in trade volume (left) and trade cost (right) for the largest trading partners. \highlight{The large relative increase in trade to India is small in absolute terms and statistically not significant.} \textbf{C} Percent change in trade volume (left bar) and change in trade cost (right bar) for selected countries.}
\label{fig:Ukraine_wheat}
\end{adjustbox}
\end{figure*}

The inference methodology presented in this work is a novel deep learning approach to solve the inverse OT problem, based on recent work on neural parameter calibration \cite{Gaskin_2023, Gaskin_2024}. We assume no underlying covariate structure \highlight{of trade costs}, but instead infer a general cost matrix $\bm{C}$, parametrized as a deep neural network, directly from data on trade flows. We train a neural network $u$ to recognise cost matrices from observations of transport plans for the global food and agricultural trade from 2000--2022 (the `\emph{training data}') by constraining it to satisfy eq.~\eqref{eq:OT}. Put simply, this means fitting the mathematical optimal transport equation to the data in such a way that the predicted cost matrices $\bm{C}(t)$ reproduce the observations $\bm{T}(t)$. The trained neural network then solves the inverse problem
\begin{equation}
    \bm{C}(t) = u(\bm{T}(t))
\end{equation}
on the observations. Though its ability to generalise to new observations depends on the amount of training data, its performance on the training data itself does not. \highlight{A probability density $\rho_C$ on the estimates is then naturally obtained as the pushforward measure
\begin{equation}
    \rho_C = u_\#\rho_T,
\end{equation}
where $\rho_T$ is the measure on $\bm{T}$. Additionally, we train a family of neural networks to capture the spread in cost matrices that optimally reproduce the transport plan  (see Methods).} As we demonstrate, this approach produces trade flow estimates that are an order of magnitude more accurate than those of a \highlight{traditional} covariate-based gravity model.

The dataset under consideration was assembled by the Food and Agricultural Organisation of the United Nations (FAO), which provides global trade matrices for over 500 products on its portal\footnote{\href{https://www.fao.org/faostat/en/\#home}{fao.org/faostat/en/\#home}} \cite{FAOStat}. Though extensive, many entries in the trade matrices are missing. Furthermore, the FAO reports two values for each bilateral flow $T_{ij}$: one reported by the exporter, and one reported by the importer. There is often a considerable discrepancy between the two, due to a multitude of epistemic factors the FAO lists in its accompanying report\footnote{\href{https://files-faostat.fao.org/production/TM/TM_e.pdf}{files-faostat.fao.org/production/TM/TM\_e.pdf}}. The uncertainty on our estimates naturally follows the uncertainty on the FAO data, without presupposing an underlying statistical model.

We apply our method to analyse global commodity flows from 2000--2022, examining the impacts of events, conflicts, trade agreements, and political changes on trade. The cost matrix uncovers economic effects that are not evident in trade volumes or retail prices alone. The article begins with a study of the war in Ukraine's impact on global wheat trade, followed by an analysis of free trade agreements and disputes in the Asia-Pacific, as well as the United Kingdom's 2016 exit from the European Union (Brexit). \highlight{We demonstrate our method's ability to provide meaningful uncertainty quantification, and compare it to a traditional gravity model, demonstrating superior prediction accuracy.}

\section{Results}
\subsection*{Case study I: the impact of the Ukrainian war on wheat trade}
The Russian Federation's invasion of Ukraine in early 2022 sent shock waves through global food markets \cite{Laber_2023}. Russia and Ukraine are two of the largest exporters of wheat, together accounting for almost 28\% of global wheat exports in 2020. The blockade of trading routes through the Black Sea and the closure or destruction of ports in Mykolaiv and Kherson meant a drop in trade to the overwhelming majority of Ukraine's export destinations, in some cases by as much as 100\% (fig.~\ref{fig:Ukraine_wheat}{\cs{\textbf{A}}}--{\cs{\textbf{B}}}). An increase of wheat exports only occurred to Europe, most significantly to \highlight{Poland, Spain, Slovakia, Romania, as well as to Algeria, and T\"urkiye}. However, our analysis shows that, although trade shrank across the globe, the accompanying increase in \highlight{wheat} trade costs disproportionately affected the Global South, in particular African nations. Of the ten countries with the largest rise in \highlight{wheat} import costs, five are in Africa, and all are in the Global South, while of the ten countries with the largest decrease in trade barriers with Ukraine, nine are in Europe (\highlight{see figs. \ref{fig:SI:Wheat_macrotrends_1}, \ref{fig:SI:Wheat_macrotrends_2}, and \ref{fig:SI:Wheat_macrotrends_3} in the SI}). Countries such as Nigeria or Angola, while experiencing a similar drop in trade as the Norway or France, simultaneously saw an increase in their trade costs. Canadian imports fell by 91\%, yet \highlight{unit trade costs} remained nearly constant, while similar drops in Zambia or the DR Congo led to marked increases in \highlight{unit trade costs, pointing to} trading barriers. European countries saw an average 9\% drop in \highlight{unit wheat trade costs} with Ukraine, while Sub-Saharan Africa saw an average 22\% increase (\highlight{see fig.~\ref{fig:SI:Wheat_macrotrends_2} in the SI}). Imports of wheat from Russia also fell globally (\highlight{see fig.~\ref{fig:SI:Russia_wheat} in the SI}), again affecting Africa particularly severely. European imports of Russian wheat fell by around 74\% with an 18\% increase in trade costs; African imports fell by on average 80\% with a 36\% increase in trade costs. While many European countries saw their imports of Ukrainian wheat rise, Russian imports fell sharply. The two largest hubs for Russian wheat, Egypt and T\"{u}rkiye, saw no change in their import volumes and small declines in their import barriers. Meanwhile, Iran saw a 0.51 point increase in \highlight{wheat unit trade costs}, leading to a 97\% percent decline in Ukrainian wheat imports. For Russian wheat, the estimated increase in trade \highlight{costs} was only 0.04, leading to a drop in imports of 53\%. Russian-Iranian trade barriers were thus not markedly affected by the war, despite a drop in trade volumes.
\begin{figure*}[th!]
\begin{adjustbox}{minipage=\textwidth-20pt}
\begin{minipage}{\textwidth}
	\cs{\textbf{A} Export of sugar products$^\dagger$ to China, selected countries}
\end{minipage}
\hfill \vspace{-2mm}

\begin{minipage}{\textwidth}
	\includegraphics[width=\textwidth]{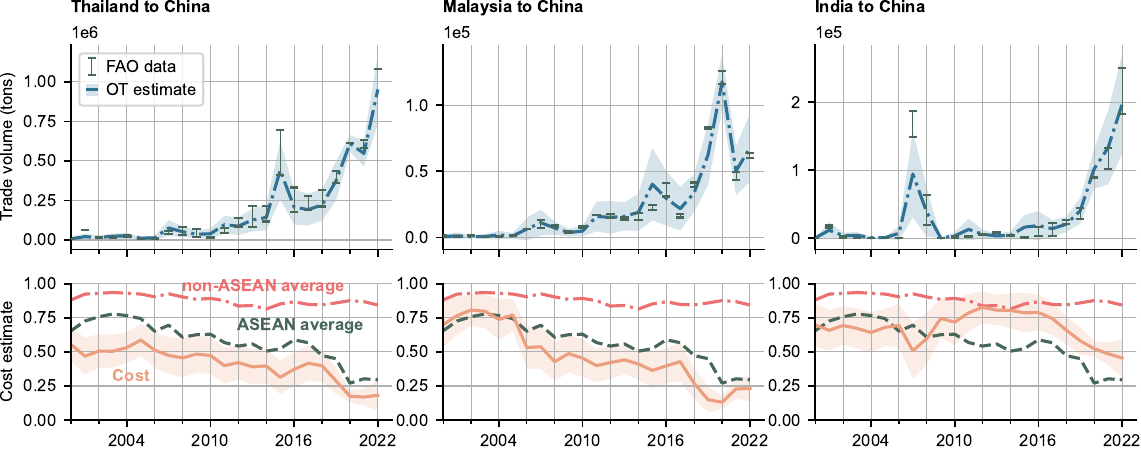}
\end{minipage}
\vspace{2mm}

\begin{minipage}{\textwidth}
	\cs{\textbf{B} Australian exports to China, selected commodities}
\end{minipage}
\hfill \vspace{-3mm}

\begin{minipage}{\textwidth}
	\includegraphics[width=\textwidth]{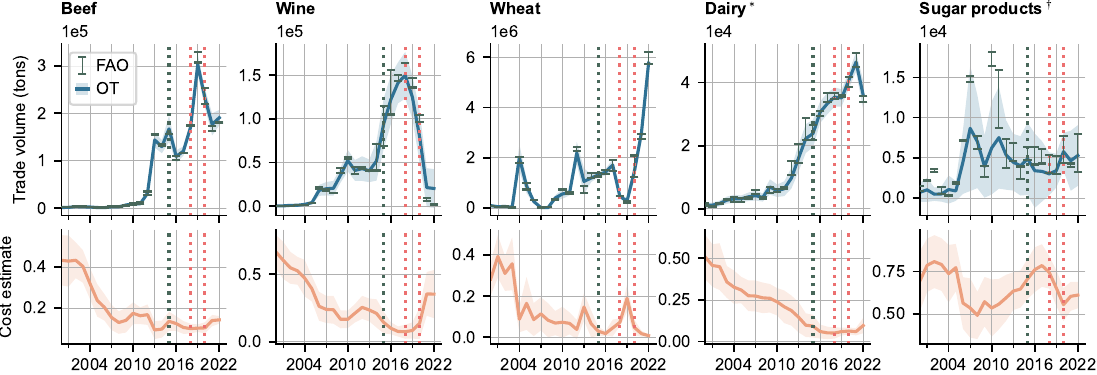}
\end{minipage}

\caption{Trade with China. \textbf{A} Export of sugar products to China. Top row: estimated trade volume (light blue) in metric tons, as well as the reported values. Bottom row: estimated cost, together with the ASEAN and non-ASEAN averages. \textbf{B} Australian exports to China, selected commodities. Top row: model estimated flow and FAO data; bottom row: estimated cost. Indicated are the signing of ChAFTA (2015, green dotted line) as well as the start of the US-China and Australia-China trade disputes (2018 and 2020, red dotted lines). Errorbands indicate one standard deviation. \smallskip \\ {\fontsize{7}{8}\selectfont $^\dagger$Sugar products comprise: sugar, refined sugar, syrups, fructose, sugar confectionery. $^\star$Dairy products comprise: butter, skim milk of cows, cheese, other dairy products.}}
\label{fig:China_trade}
\end{adjustbox}
\end{figure*}
\begin{figure*}[ht!]
\begin{adjustbox}{minipage=\textwidth-20pt}
\begin{minipage}{\textwidth}
	\cs{\textbf{A} Global barley trade, 2018 (left) and 2021}
\end{minipage}
\hfill \vspace{-3mm}

\begin{minipage}{\textwidth}
	\includegraphics[width=\textwidth]{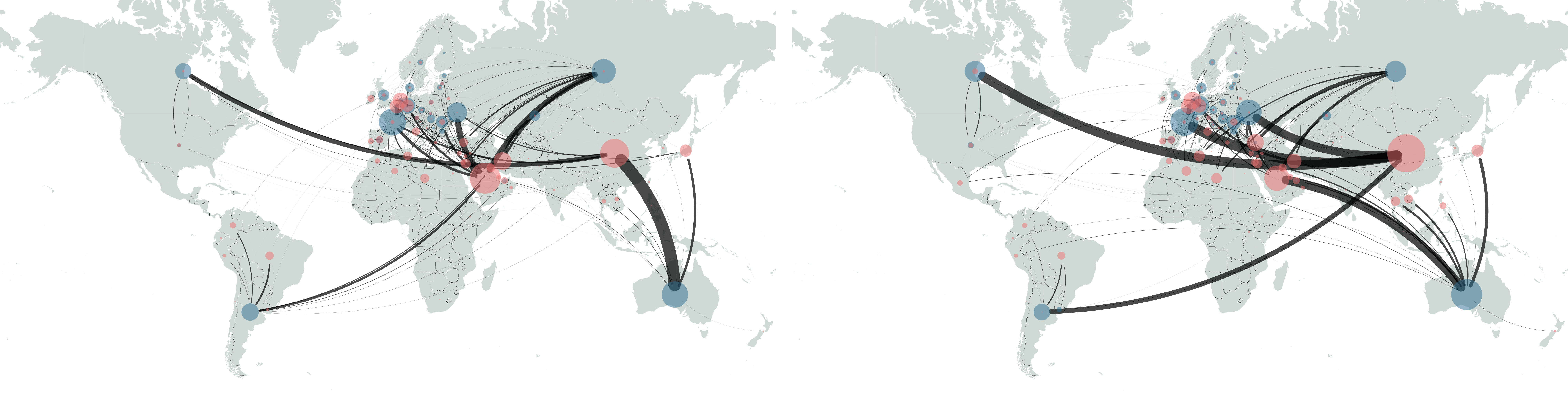}
\end{minipage}

\begin{minipage}{\textwidth}
	\cs{\textbf{B} Barley export cost to China, 2018 (left) and 2021}
\end{minipage}
\hfill \vspace{-3mm}

\begin{minipage}{\textwidth}
	\includegraphics[width=\textwidth]{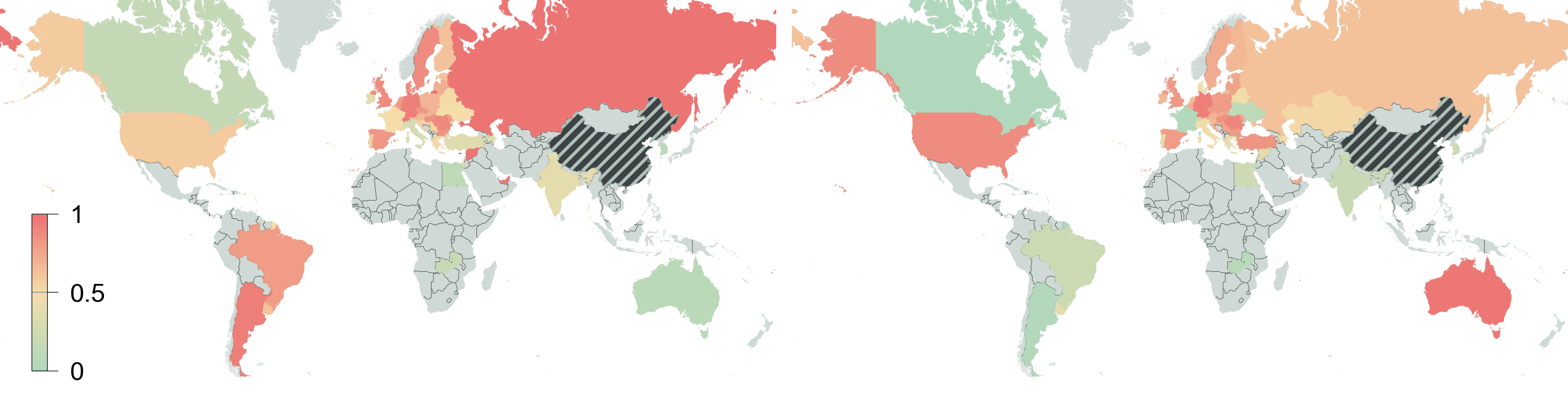}
\end{minipage}
\vspace{2mm}

\begin{minipage}{\textwidth}
	\flushleft \cs{\textbf{C} Trade volumes and costs}\vspace{1mm}
\end{minipage}
\begin{minipage}{\textwidth}
	\includegraphics[width=\textwidth]{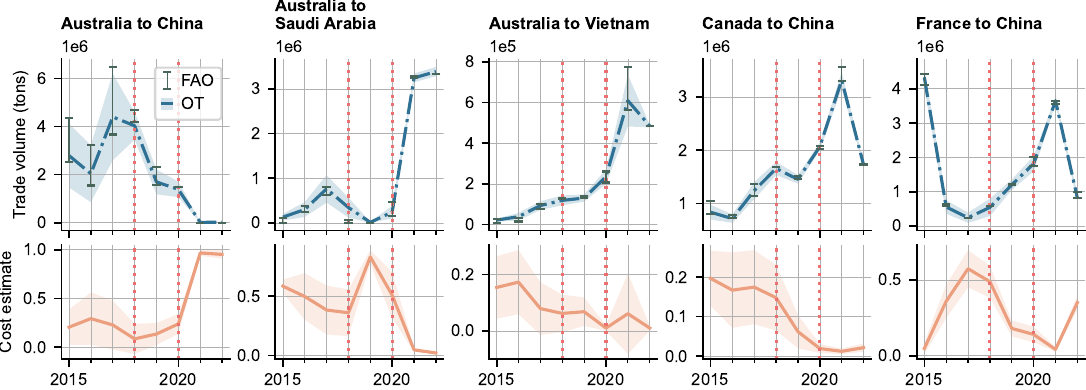}
\end{minipage}

\caption{Global barley trade between 2015--2022. After the introduction of Chinese import tariffs on Australian barley in May 2020, the entire supply chain restructured itself, with Chinese demand being supplied from France, Canada, and Ukraine, and Australia increasingly exporting to Saudi Arabia and Southeast Asia. \textbf{A}--\textbf{B} Trade in in metric tons, 2018 and 2021. Import values are shown in red, export values in blue. \textbf{C} Model estimated trade volumes (top row) and cost (bottom row) for selected countries. Dotted lines indicate the start of the US-China and Australia-China trade wars.}
\label{fig:Barley_trade}
\end{adjustbox}
\end{figure*}
\begin{figure*}[ht!]
\begin{adjustbox}{minipage=\textwidth-20pt}
\begin{minipage}{\textwidth}
	\cs{\textbf{A} Global soya bean trade, 2016 (left) and 2018}
\end{minipage}
\hfill \vspace{-3mm}

\begin{minipage}{\textwidth}
	\includegraphics[width=\textwidth]{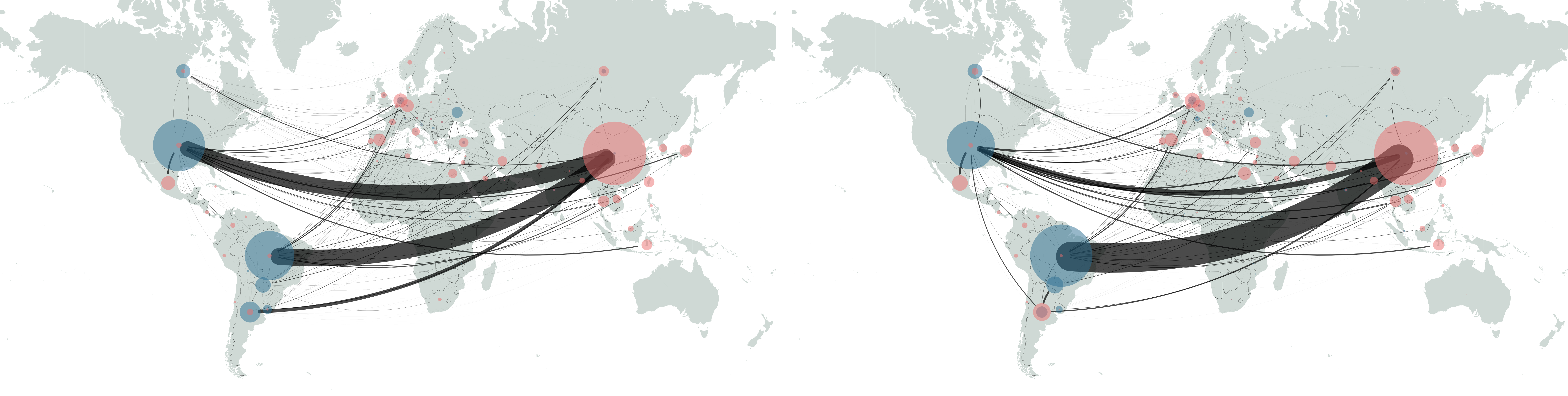}
\end{minipage}

\begin{minipage}{\textwidth}
	\flushleft \cs{\textbf{B} Yield}
\end{minipage}
\hfill \vspace{-3mm}

\begin{minipage}{\textwidth}
	\includegraphics[width=\textwidth]{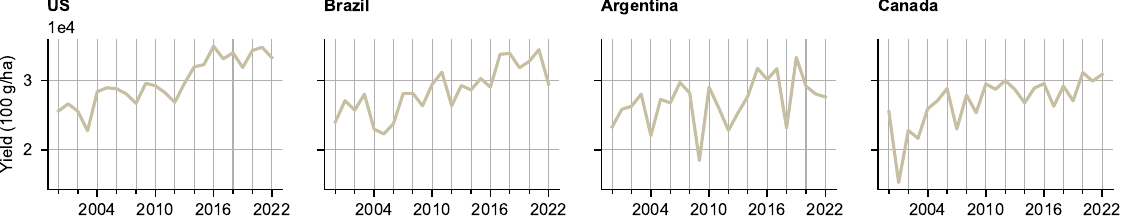}
\end{minipage}
\vspace{1mm}

\begin{minipage}{\textwidth}
	\cs{\textbf{C} Trade volumes and utilities}
\end{minipage}
\hfill \vspace{-3mm}

\begin{minipage}{\textwidth}
	\includegraphics[width=\textwidth]{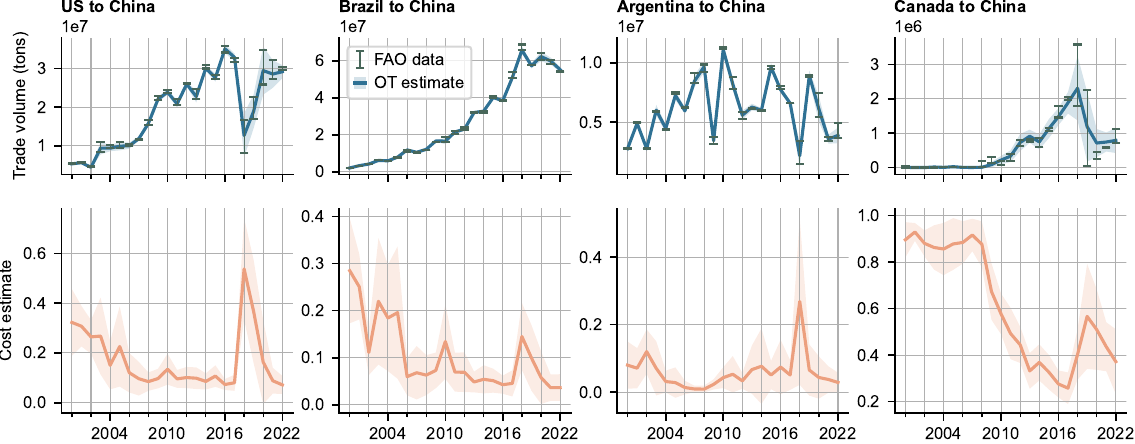}
\end{minipage}
\vspace{1mm}

\caption{Global soya bean trade. \textbf{A} In 2018, the Chinese government raised import tariffs on American soya beans in a retaliatory action against US trade restrictions. The shortfall was met by imports from Brazil. \textbf{B} Soya bean yield in 100 g/hectare. Argentina in 2018 experienced a major drop in yields, leading to an increase in exports from the US. \textbf{C} Predicted trade volumes in metric tons (top row) and predicted cost (bottom row).}
\label{fig:Soya_trade}
\end{adjustbox}
\end{figure*}
\begin{figure*}[ht!]
\begin{adjustbox}{minipage=\textwidth-20pt}
\begin{minipage}{\textwidth}
	\cs{\textbf{A} Vegetables}
\end{minipage}
\begin{minipage}{\textwidth}
	\includegraphics[width=\textwidth]{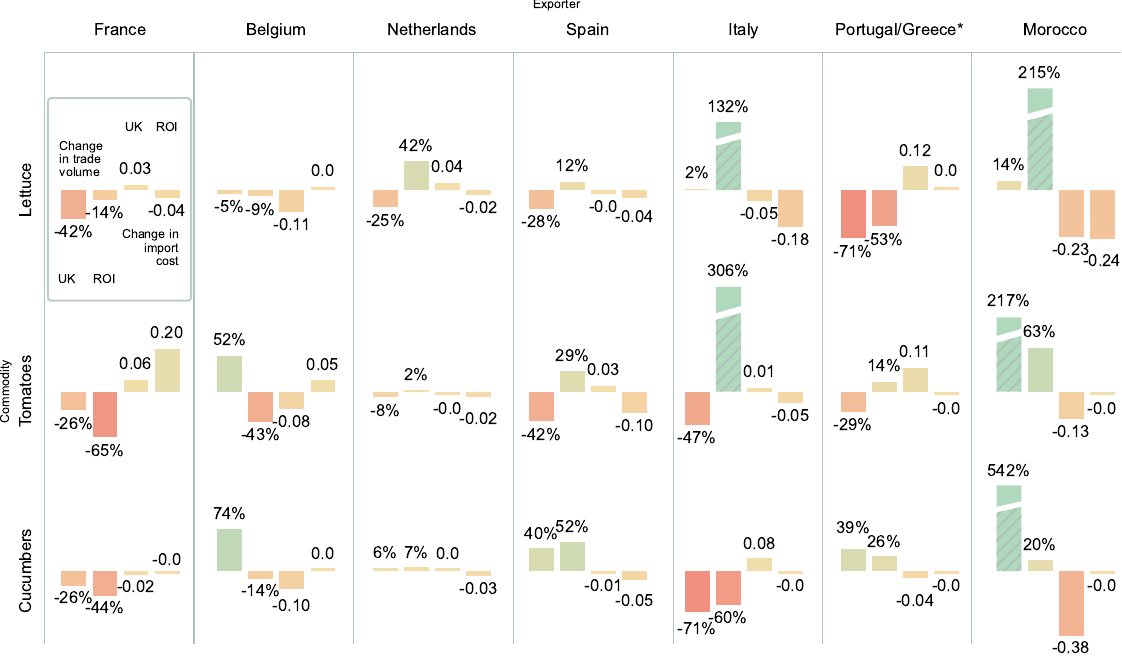}
\end{minipage}
\vspace{1mm}

\begin{minipage}{\textwidth}
	\cs{\textbf{B} Wine} 
\end{minipage}
\begin{minipage}{\textwidth}
	\includegraphics[width=\textwidth]{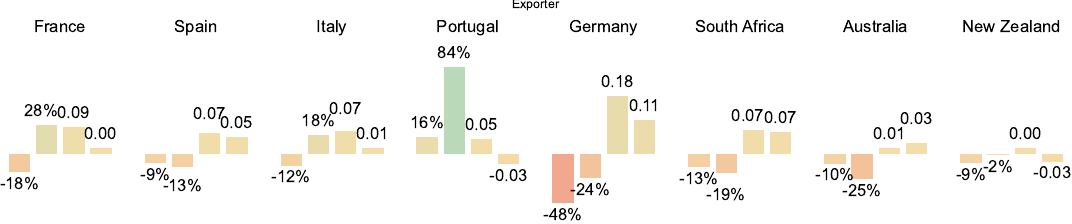}
\end{minipage}
\caption{Change in UK and Ireland (ROI) imports, 2016--2022. For each exporting country, the left two bars indicate the percent change in trade volume between 2016 and 2022 for the UK and the ROI respectively, the right two bars show the change in import costs. \textbf{A} Vegetable imports. Top row: lettuce (including chicory), and other fresh vegetables; middle row: tomatoes; bottom row: cucumbers and gherkins. \textbf{B} Wine imports. \smallskip \\ {\fontsize{7}{8}\selectfont $^*$Exporter for cucumbers is Greece.}}
\label{fig:Brexit}
\end{adjustbox}
\end{figure*}
\subsection{Case study II: Trade in Southeast Asia and Asia-Pacific}
A series of free-trade agreements came into effect in Southeast Asia and the Asia-Pacific region in the 2000s and 2010s, significantly among them the China-Australia Free Trade Agreement (ChAFTA) in 2015, the ASEAN-China free trade agreement (ACFTA, gradually entering into force from 2003) and the Comprehensive and Progressive Agreement for Trans-Pacific Partnership (CPTPP) between 11 counties bordering the Pacific Ocean (2018) \cite{ChAFTA, ACFTA, CPTPP}. Together with China's accession to the WTO in 2001 and its rapid economic growth, these trade agreements coincide with some of the largest increases in trade flows in recent history. In figure \ref{fig:China_trade}{\cs{\textbf{A}}}, we show the trade flow of sugar and sugar products from Thailand, Malaysia, and India to China, as well as the estimated costs. In our model, the cost of importing sugar from Thailand fell consistently from 2000--2022, following a general trend for ASEAN countries (bottom row, green line) which commenced around 2005.  Indian exports, by comparison, remained relatively low until 2015, when Indian prime minister Narendra Modi visited China, and top officials from both sides agreed to increase bilateral trade to US\$100 billion by the end of the year. This visit marked a dramatic shift in Indo-Chinese trade, as exemplified by the huge increase in sugar trade. From 2015--2022, \highlight{sugar} export cost from India dropped sharply by 33\%, precipitating a steep increase in trade. By contrast, trade cost from non-ASEAN members has remained constant over the past twenty years (red line, fig.~\ref{fig:China_trade}{\cs{\textbf{A}}}).

The PRC is one of Australia's largest export markets for food and agricultural products. Our analysis shows a precipitous reduction in trade barriers for Australian exports since China's accession to the WTO in 2001 (see fig.~\ref{fig:China_trade}{\cs{\textbf{B}}}), particularly for beef, wheat, wine, and dairy. Between 2002 and 2010, these commodities saw a 30--50\% drop in \highlight{their respective trade costs}. Our estimates indicate that ChAFTA had little effect on Australian trade costs, since it succeeded a period of deepening ties. Dairy \highlight{trade costs}, for instance, had already fallen from 0.51 to 0.1 from 2000 to 2015, thereafter only falling a further 0.04 points until 2020. Wine exports too saw their largest reductions in \highlight{trade costs} between 2000 and 2010, only experiencing a 0.07 drop from 2015 to 2018 compared to the 0.52 point reduction from 2000--2015.

In January 2018, the \highlight{first} Trump administration started imposing import tariffs on goods primarily from China. In response, the Chinese government increased tariffs on a variety of products, including agricultural imports. The largest agricultural export from the US to China, soya beans, were hit with a 25\% import tariff \cite{Biesheuvel_2018}. Meanwhile, political tensions between China and Australia caused Beijing to introduce high anti-dumping tariffs on Australian exports such as barley (80.5\%) and wine (206\%), starting in 2020 \cite{Sullivan_2022}. Wine trade had previously been tariff-free since the signing of ChAFTA in 2015. Our analysis provides an estimate of the change in the ease of trading these measures induced (figs. \ref{fig:China_trade}{\cs{\textbf{B}}} and \ref{fig:Barley_trade}). Australian beef, wine and barley imports all experienced large increases in cost, following the implosion of trade volumes. Australia was able to divert some of its excess barley supply to Saudi Arabia, which saw a decrease in \highlight{trade costs} of over 0.8 points between 2019 and 2022 (fig.~\ref{fig:Barley_trade}{\cs{\textbf{C}}}). Trade volumes to Vietnam also increased from 200,000 to 800,000 metric tons, though trade costs remained approximately constant. Meanwhile, after 2020 China doubled its barley imports from Canada and France. We found that import barriers from both countries were reduced only slightly in 2021 and rebounded the following year.
\begin{figure*}[ht!]
\begin{adjustbox}{minipage=\textwidth-20pt}
\begin{minipage}{0.5\textwidth}
    \cs{\textbf{A} Estimates vs FAO data, cucumbers} \vspace{2mm}
\end{minipage}
\begin{minipage}{0.5\textwidth}
    \cs{\textbf{B} Estimates vs FAO data, barley} \vspace{2mm} 
\end{minipage} 
\begin{minipage}{0.5\textwidth}
	\includegraphics[width=\textwidth]{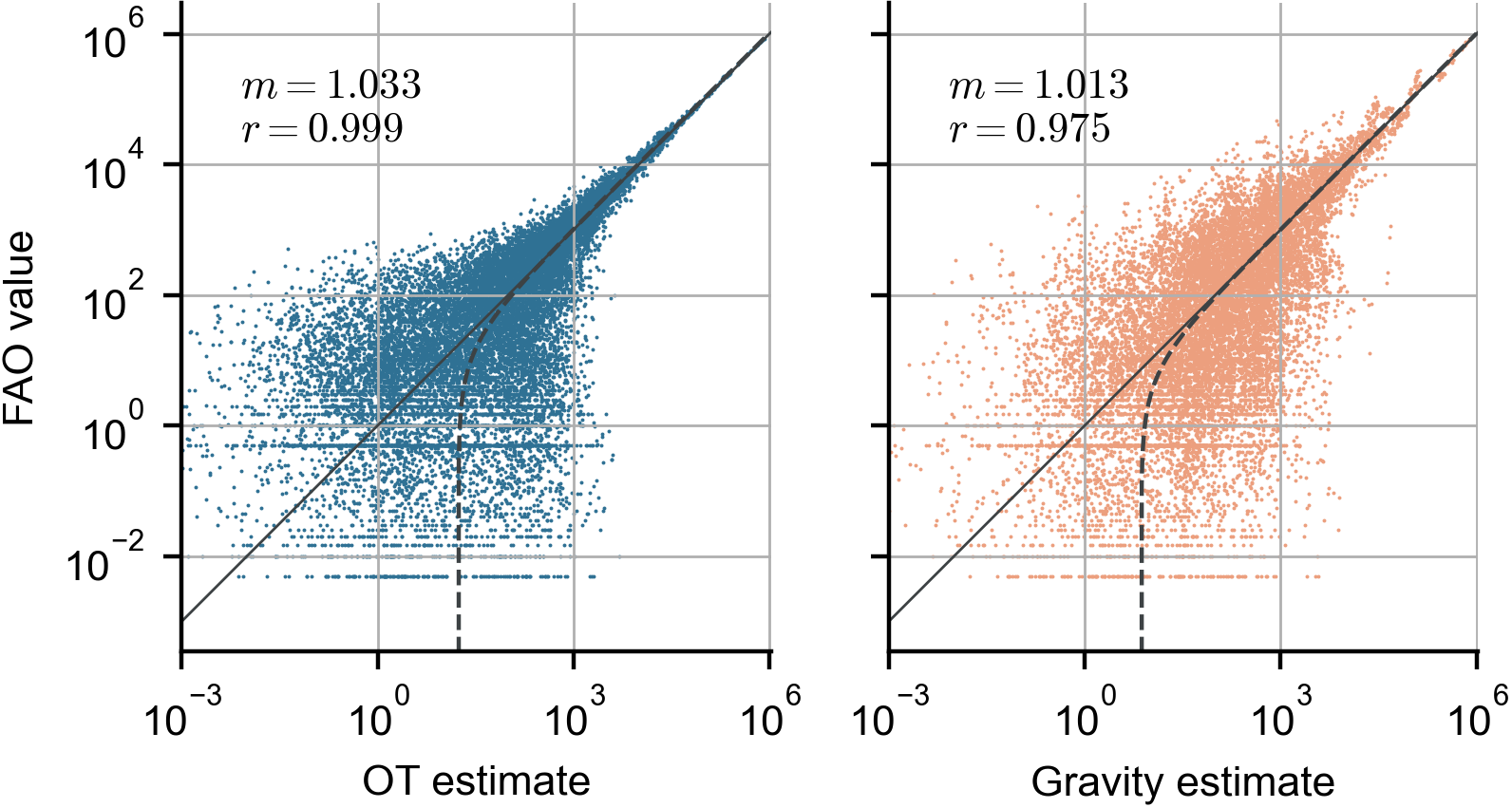}
\end{minipage}
\begin{minipage}{0.5\textwidth}
	\includegraphics[width=\textwidth]{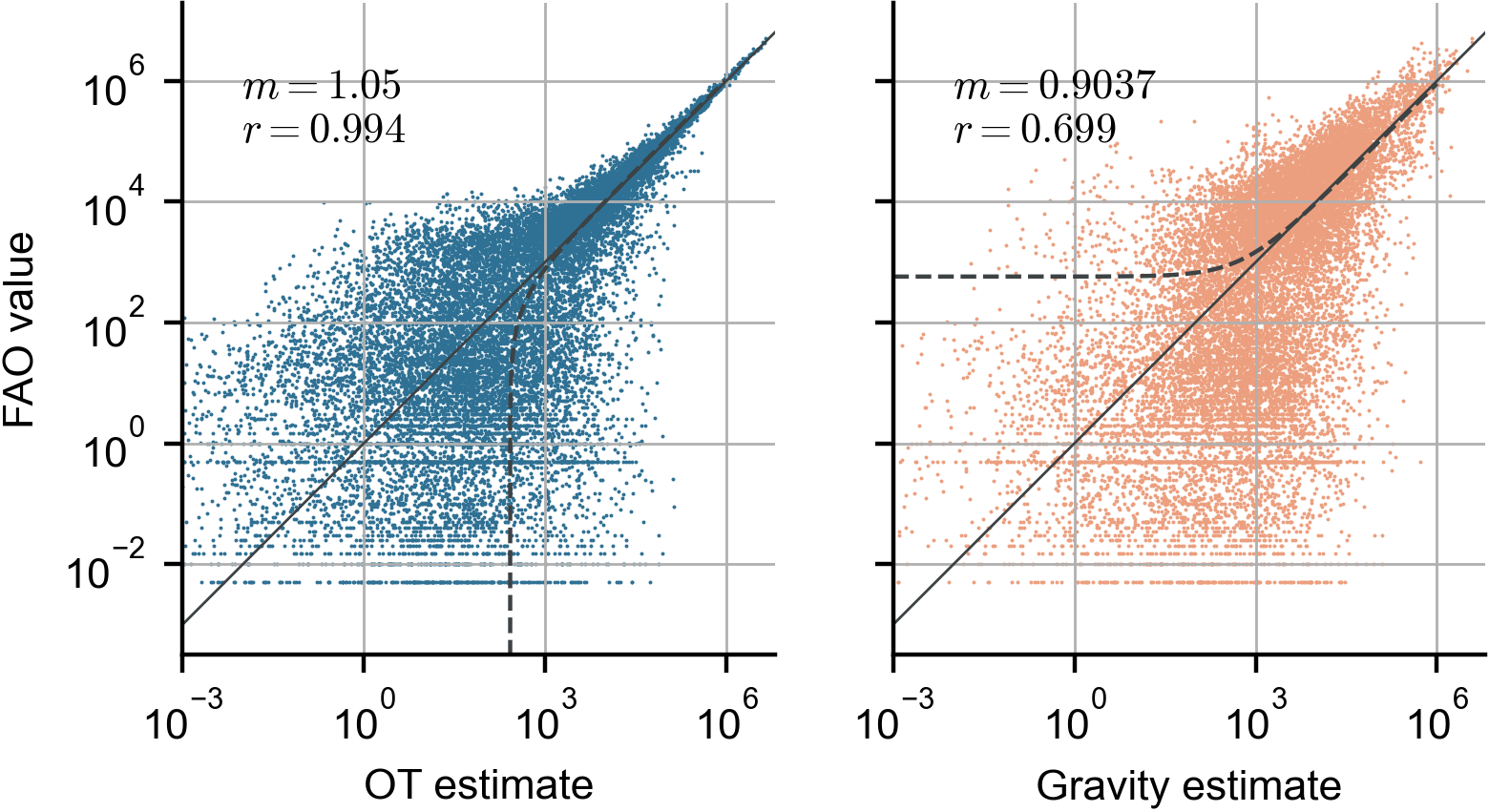}
\end{minipage}

\vspace{4mm}
\begin{minipage}[ht]{0.46\textwidth}
	\cs{\textbf{C} Estimation accuracy (RMSE)}
\end{minipage}
\begin{minipage}[ht]{0.54\textwidth}
	\cs{\textbf{D} Estimation accuracy (RMSE in standard deviations)}
\end{minipage}
\hfill \vspace{-2mm}

\begin{minipage}{\textwidth}
	\includegraphics[width=\textwidth]{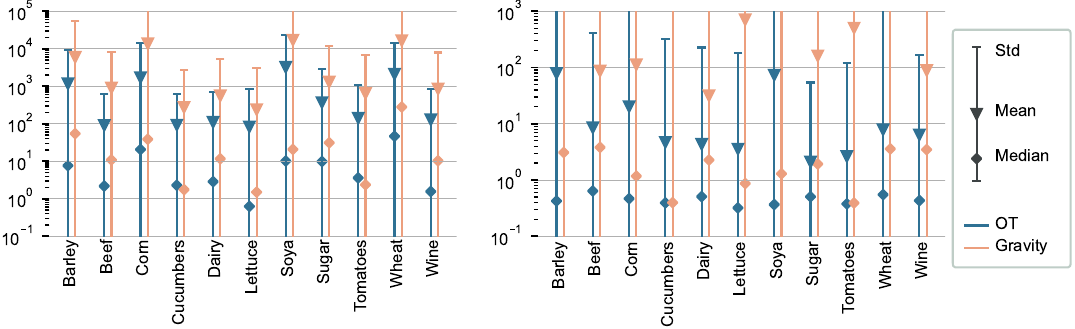}
\end{minipage}
\caption{Comparison with gravity model. \textbf{A}--\textbf{B} Comparison plot of the OT and gravity estimates ($x$-axis) versus the true data ($y$-axis) on two selected commodities. Also shown is a linear fit (dotted line), its estimated slope $m$, the Pearson coefficient of the fit $r$, and the line $y=x$ (solid line). See fig.~\ref{fig:SI:Gravity_comparisons_all} in the appendix for an overview of all commodities. \textbf{C}--\textbf{D} Comparison of the RMSE accuracies of the estimated transport volumes of the OT approach (blue) and the gravity model (orange). Values are averaged over all countries and years, with the errorbars showing one standard deviation from the mean (triangular marker). Also shown are the median values (diamond markers). Shown are the RMSE ({\cs{\textbf{C}}}) and the RMSE in units of the standard deviation on the true data ({\cs{\textbf{D}}}).}
\label{fig:Gravity_comp}
\end{adjustbox}
\end{figure*}

\subsection{Case study III: Brexit}
In 2016, the United Kingdom voted to leave the European Union, officially exiting the common market and customs union on December 31, 2020. This case study examines the impact of Brexit on British import patterns by comparing vegetable and wine imports from mainland Europe to both the United Kingdom and the Republic of Ireland (ROI), which remains part of the Eurozone and the common market. While both island nations naturally source the majority of their fresh produce from mainland Europe, their trading patterns have evolved in markedly different ways (fig.~\ref{fig:Brexit}{\cs{\textbf{A}}}). Imports of lettuce from Europe generally fell for the UK, accompanied by a rise in import cost: 25\% decrease in trade volume and +0.04 in import costs from the Netherlands, the largest exporter of lettuce and chicory to the UK, as well as a –28\% drop in trade from Spain, though with no change in import cost. Ireland increased its imports of lettuce and other greens from the Netherlands, Spain, and Italy, accompanied by a general decrease in trading costs \highlight{for those products}. \highlight{The ROI's imports of greenery from Portugal fell by 53\%, with no change in trade barriers; for the UK, a 71\% fall in trade volumes accompanied a 0.12 point increase in trade cost}. In the case of the Netherlands, Ireland saw a consistent reduction in vegetable trade costs, unlike the United Kingdom. It is interesting to note that the United Kingdom significantly increased its imports of vegetables from Morocco, accompanied by a precipitous drop in trade costs, indicating a facilitation of trade between the two countries in the wake of Brexit. This is not true for the ROI: though it increased its imports of Moroccan tomatoes and cucumbers, Irish trade costs remained mostly unchanged.

A more clear-cut trend emerges in the wine trade (fig.~\ref{fig:Brexit}{\cs{\textbf{B}}}): here, the UK was consistently affected more negatively than the Republic of Ireland: British import costs from all eight countries considered rose by considerably more than those of the ROI. A 9\% drop in Spanish wine import \highlight{by the UK} was accompanied by a 7\% increase in trading costs, while a \highlight{13\% drop} in Irish imports only meant a 5\% increase in costs. Portugese wine imports to the UK rose by 16\%, notwithstanding a 0.05 increase in trade costs. A similar pattern holds for South African, Australian, and New Zealand imports. The EU maintains free trade or regulatory agreements removing wine import duties with the former two \cite{EU-SADC, WFA_2016}. When the UK left the European Union, wine from Australia entered the UK Global Tariff rate, which in mid-2023 was eliminated under the Australia-United Kingdom FTA \cite{DAFF_2020}. South African wines, by contrast, continued to be imported to the UK tariff-free post-Brexit \cite{VINPRO_2020}. Yet here too, the United Kingdom's 13\% decrease in imports was driven by a 0.07 increase in trading costs; Ireland, by contrast, imported 19\% less wine, driven by a comparable 0.07 point increase in trading costs.

\subsection{Comparison with Gravity model}

Lastly, we compare the performance of our method with a \highlight{traditional} gravity model \cite{Yotov_2017, CEPII, UNCTAD_TRAINS, WTO_Stats}, as specified in equations \eqref{eq:Gravity}--\eqref{eq:Gravity_covariates}. The covariates include geographic distance, shared borders, colonial ties, common language, regional trade agreements, tariffs, and importer/exporter fixed effects. We estimate the coefficients using Poisson Pseudo Maximum Likelihood estimation \cite{silva2006log,Correia2020} and compare the accuracy of the estimated transport plans $\bm{T}$. Figures \ref{fig:Gravity_comp}{\cs{\textbf{A}}}--{\cs{\textbf{B}}} show scatter plots of the OT (blue) and the gravity (orange) estimates against the FAO data. For all commodities studied, a linear fit through the OT estimates yields a near-perfect slope of $m=1$ with a Pearson coefficient close to $1$, \highlight{particularly providing an exact fit for the upper tail of the trade value distribution}. In contrast, the gravity model's performance is much more volatile, with linear fits ranging from a Pearson coefficient of between 0.975 (best) to 0.699 (worst) (see also SI). Due to model misspecification, the fits to the tails of the distributions are generally significantly poorer. Consequently, figure \ref{fig:Gravity_comp}{\cs{\textbf{C}}} shows that the OT approach significantly outperforms the gravity model in terms of RMSE, \highlight{often by an order of magnitude}. Figure \ref{fig:Gravity_comp}{\cs{\textbf{D}}} illustrates that OT estimates typically fall within one standard deviation of the data uncertainty, whereas gravity estimates tend to range from one to two, at times even three to four standard deviations. The gravity model also exhibits much higher variance in accuracy compared to OT. \highlight{We further an investigate an alternative specification by replacing time-varying country-level and time-invariant pair regressors with exporter-time, importer-time, and exporter-importer fixed effects in a three-way gravity model \cite{Breinlich2021}, capturing multilateral resistance terms more effectively (see eq.~\eqref{eq:Gravity_model_alt} in the SI for details). As this formulation incorporates higher degrees of freedom in the three-way specification, its performance naturally converges toward that of OT, which more accurately fits high-value trade flows. Incorporating time-varying bilateral trade regressors would likely further narrow this gap, reinforcing the validity of the OT approach. Notably, optimal transport provides trade cost estimates even without such regressors---many of which are difficult to quantify.}
\newpage
\section{Discussion}
This paper introduces a novel and versatile approach for identifying the drivers and barriers of global commodity trades. Using optimal transport theory, we are able to obtain a cost structure that is more expressive than a covariate-based gravity approach. Our estimates are thus orders of magnitude more accurate than \highlight{traditional gravity models, matching the performance of high-dimensional fixed-effects specifications while maintaining consistent accuracy across datasets}. The optimal transport approach models trade networks as a dynamical, interconnected system, allowing to capture complex rearrangements and network response dynamics to e.g. trade wars, conflicts, or shifts in political relations. Though the current work looks only at global agrifood markets, the methodology proposed is general and applicable to commodity flows, financial markets, or banking networks \cite{Giesecke2019}. Beyond economics, the optimal transport approach also relates e.g. to global migration flows, which can be estimated from migrant stock data \cite{Willekens_1999,Azose_2018,Gaskin_Abel_2025}. \highlight{Future work could explore correlations between related commodities within this framework, as well as develop hybrid models that combine observed covariates with data-driven residuals via semi-structured OT costs. Another promising direction is counterfactual analysis using the conditional equilibrium framework of \cite{Yotov_2017}. Extending the model to a full general equilibrium setting with endogenous production and pricing is also a natural next step.}

\section{Methods}
\subsection{Entropy-regularised Optimal Transport}

In OT one wishes to find the optimal flow of mass from a source distribution to a target distribution, while minimising an overall transport cost. This abstract problem has a wide range of applications in economics, logistics, image restoration, transport systems, or urban structure \cite{Galichon_2016, Santambrogio_2015, Peyre_2019}.

Consider an $m$-dimensional space $X$, an $n$-dimensional space $Y$, and $\bm{C}$ a measure on $X \times Y$. The entries of  $\bm{C}$ correspond to the cost of transporting mass from one location in $X$ to a target in $Y$. Given two probability measures $\bm{\mu} \in P(X)$ and $\bm{\nu} \in P(Y)$ (the supply and demand), the OT problem consists in finding a \emph{transport plan} $\bm{T}$ minimising the overall cost eq.~\eqref{eq:OT_cost}. The transport plan $\bm{T}$ must also satisfy the \emph{marginal constraints}
\begin{align}
    \sum_j T_{ij} = \bm{\mu} \text{ and } \sum_i T_{ij} = \bm{\nu}.
    \label{eq:marginals}
\end{align}
In practice one usually solves the \emph{entropy regularised} OT formulation, which can be solved much more efficiently \cite{Cuturi_2013}; here, an additional term is added to the objective:
\begin{align}
\min_{\bm{T}} \sum_{ij} C_{ij} T_{ij} + \varepsilon \sum_{ij} T_{ij} \left( \log T_{ij} - 1 \right),
\label{eq:OT_entropy_regularised}
\end{align}
where $\varepsilon>0$ is a positive regularisation parameter. This regularisation prevents monopolisation, i.e. demand being supplied from only a few sources.

This constrained optimisation problem eq.~\eqref{eq:OT_entropy_regularised} can be solved by considering the Lagrangian
\begin{align}
    \mathcal{L} &= \sum_{ij} T_{ij}C_{ij} + \langle \bm{\lambda}, \sum_j T_{ij} - \bm{\mu} \rangle + \langle \bm{\eta}, \sum_i T_{ij} - \bm{\nu} \rangle \nonumber \\
    & + \varepsilon \sum_{ij} T_{ij} \left( \log T_{ij} - 1 \right),
\end{align}
with $\bm{\lambda} \in \mathbb{R}^m$ and $\bm{\eta} \in \mathbb{R}^n$ Lagrangian multipliers. Minimising $\mathcal{L}$ with respect to $\bm{T}$ gives the solution
\begin{equation}
    T_{ij} = e^{-\lambda_i / \varepsilon} e^{-C_{ij}/ \varepsilon} e^{-\eta_j /\varepsilon}
\end{equation}
or 
\begin{equation}
    \bm{T} = \bm{\Pi} e^{-\bm{C} / \varepsilon} \bm{\Omega},
    \label{eq:OT_solution}
\end{equation}
where $\bm{\Pi} = \mathrm{diag}(e^{-\lambda_1 / \varepsilon}, \cdots, e^{-\lambda_m / \varepsilon}) \in \mathbb{R}^{m \times m}$ and $\bm{\Omega} = \mathrm{diag} (e^{-\eta_1 / \varepsilon}, \cdots, e^{-\eta_n / \varepsilon}) \in \mathbb{R}^{n \times n}$ are diagonal matrices of Lagrangian multipliers. 

Finding $\bm{\Pi}$ and $\bm{\Omega}$ is achieved through an iterative scaling procedure that is variously called \emph{Iterative Proportional Fitting} (IPF), \emph{RAS}, or \emph{Sinkhorn's algorithm} \cite{Deming_1940, Cuturi_2013, Idel_2016}. Define $\bm{M} = e^{-\bm{C}/\varepsilon}$; then, given an initial guess $\bm{\Pi}^0$, we update $\bm{\Omega}$ to satisfy the first marginal constraint eq.~\eqref{eq:marginals} 
\begin{equation}
    \bm{\Omega} \bm{M}^\top \bm{\Pi}^0 = \bm{\nu}.
\end{equation}
Solving for $\bm{\Omega}$ gives
\begin{equation}
    \bm{\Omega}^0 = \dfrac{\bm{\nu}}{\bm{M}^\top \bm{\Pi}^0}
\end{equation}
where the division is understood element-wise. Similarly, we obtain the next update for $\bm{\Pi}$ as
\begin{equation}
    \bm{\Pi}^1 = \dfrac{\bm{\mu}}{\bm{M} \bm{\Omega}^0},
\end{equation}
and so on. 
The algorithm can thus be summarised as follows:
\begin{algorithm}[H]
\caption{Sinkhorn's Algorithm}
\label{alg:Sinkhorn}
\begin{algorithmic}[1]
    \Inputs{$\bm{M}$ (Exponential of cost matrix) 
    \\ $\bm{\mu}, \bm{\nu}$ (marginals) 
    }
    \State {Initialise the first Lagrangian multiplier $\bm{\Pi}^0$}
    \For {$n$ iterations}
    \State {$\bm{\Omega}^{i+1} \gets \dfrac{\bm{\nu}}{\bm{M}\bm{\Pi}^i}$}
    \State {$\bm{\Pi}^{i+1} \gets \dfrac{\bm{\mu}}{\bm{M}^\top\bm{\Omega}^{i+1}}$}
    \EndFor
\end{algorithmic}
\end{algorithm} 
\noindent Under certain conditions, convergence of the algorithm to a unique solution is guaranteed \cite{Sinkhorn_1964, Rueschendorf_1995}. 

\highlight{The classic OT problem eqs.~\eqref{eq:OT_cost}--\eqref{eq:OT_constraints} can be interpreted as the central planner's problem of finding the optimal assignment or matching of supply and demand. The entropy-regularised OT problem can be viewed as a similar optimal assignment problem, but subject to uncertainty and/or randomisation. The dual problem to entropy-regularised OT, i.e. minimising \eqref{eq:cost_eps} subject to \eqref{eq:OT_cost}, is given by
\begin{align}
    \label{eq:dualOT}
    \max_{\bm{f}, \bm{g}} \langle \bm{f}, \bm{\mu}\rangle + \langle \bm{g}, \bm{\nu}\rangle\highlight{ - \varepsilon \sum_{i,j} e^{\frac{f_j + g_i - C_{ij}}{\varepsilon}}}.
\end{align}}
\highlight{In the limit $\varepsilon \rightarrow 0$ the last term ensures that the dual potentials $\bm{f}$ and $\bm{g}$ satisfy }
\begin{align}\label{eq:dual_constraint}
\bm{f} \oplus \bm{g} \leq \bm{C}
\end{align}
where $\bm{f} \oplus \bm{g} = \bm{f} \bm{1}^T_m + \bm{g} \bm{1}^T_n$. \highlight{Condition \eqref{eq:dual_constraint} corresponds to the admissibility condition of the dual non-regularised OT problem. In this context} $\bm{f}$ and $\bm{g}$ can be interpreted as the minimal cost of picking up and dropping off a good at locations respectively. The problem of finding the optimal plan $\bm{T}$ is thus split into determining the optimal cost of collecting and delivering goods. The constraint \eqref{eq:dual_constraint} ensures optimality: if $f_i + g_j > C_{ij}$, that is, the cost of picking up a good at location $i$ and dropping it off at location $j$ is larger than the transportation cost, it cannot be optimal.

\subsection{Neural Inverse Optimal Transport}
\highlight{In inverse optimal transport one wishes to infer the underlying cost $\bm{C}$ from (partial) observations of transport plans $\bm{T}$, which are usually solutions to entropy regularised OT problems. Rewriting eq.~\eqref{eq:OT_solution}, we have
\begin{equation*}
    \bm{C} = \varepsilon \left( \log \bm{\Pi} + \log \bm{\Omega}^\top - \log \bm{T} \right);
\end{equation*}
thus $\bm{C}$ is only determined up to an additive decomposition into row and column vectors, with any transformation of the kind
\begin{equation*}
    C_{ij} \mapsto C_{ij} + \alpha_i + \beta_j
\end{equation*}
leaving the transport plan $\bm{T}$ invariant, since the transformation can be absorbed by the scaling vectors. To constrain the problem, we can bound the cost $C_{ij} \in [0, C_{max}]$, and demand that the inverse of $0$ should be $C_{max}$, i.e.
\begin{equation*}
    u: \bm{T} \mapsto \bm{C}, \, u(0) = C_{max}.
\end{equation*}
This is a natural restriction, since it implies that the cost on edges with zero transport flow should be maximal. If this maximum is attained in every row and every column of the cost matrix, the row- and column-shifts $\alpha_i$ and $\beta_j$ must satisfy
\begin{equation*}
    C_{max} = \max_j C_{ij} \overset{!}{=} \max_j (C_{ij} + \alpha_i) \leq C_{max} \ \forall \ i \Rightarrow \alpha_i =0.
\end{equation*}
(and similarly $\beta_j$). Thus, under these conditions the cost matrix is uniquely determined (see also fig.~\ref{fig:SI:validation} in the SI).}

To infer the cost matrix function $\bm{C}(t)$ from a dataset of transport plan observations $\bm{T}(t)$, we build on the neural parameter estimation method first introduced in \cite{Gaskin_2023} and subsequently expanded upon \cite{Gaskin_2024}. We wish to train a neural network $u$ to solve the inverse OT problem $\bm{C}(t) = u(\bm{T}(t))$. We do so by constructing a loss function \highlight{based on the optimal transport equations, i.e.
\begin{equation}
	J = \Vert \predm{T}(\predm{C}) - \bm{T} \Vert_2^2 + \sum_{(i, j) \in \mathcal{S}} (C_{ij} - C_{max})^2.
	\label{eq:loss_function}
\end{equation}}
Here, $\predm{T}(\predm{C})$ is the estimated transport plan obtained by solving Sinkhorn's algorithm alg. [\ref{alg:Sinkhorn}] \highlight{until convergence (determined by a numerical tolerance criterion), and
\begin{equation*}
    \mathcal{S} := \{(i, j) \ \vert \ T_{ij} = 0\}
\end{equation*}
are the zero-flow edges of $\bm{T}$. The second term thus enforces the maximum value $C_{max}$ to be attained on $\mathcal{S}$. Crucially, the solution of entropy-regularised OT is differentiable with respect to its inputs, and the derivative of $\predm{T}(\predm{C})$ with respect to $\predm{C}$ can be computed numerically. Thus the loss $J$ can be minimised using gradient descent methods.} The data is processed in batches, and a gradient descent step performed on the neural network parameters after each batch. The loss is only calculated for links with trade flow $>0$. \highlight{Note that our goal is not to predict future trade-flows, but rather to infer an underlying cost which drives the flows subject to the optimal transport model. We therefore do not require large volumes of data, as would be typical in a prediction task.}  

As mentioned, the FAO dataset contains two values for each entry $T_{ij}$: one reported by the exporter, and one by the importer. Let $\bm{T}^\text{E}$ be the transport plan where all entries are those reported by the exporters, and $\bm{T}^\text{I}$ those where all are reported by the importers. The \emph{training data}---i.e., the data we use to train the function $u$---consists of only these two transport plans for each year: $\{\bm{T}^\text{E}(t), \bm{T}^\text{I}(t)\}$, giving a total training set size of $2 \times L$, where $L=23$ are the number of observation points. A hyperparameter sweep showed that using a deep neural network with 5 layers, 60 nodes per layer, and hyperbolic tangent activation functions on all layers but the last, where we use a sigmoid, gives best results. Using a sigmoid activation function on the last layer ensures $0 \leq C_{ij} \leq 1 = C_{max}$. We use the Adam optimizer \cite{Kingma_2014} to train the neural network. We pool all FAO trade matrices to only contain those countries that account for 99\% of import and export volumes, subsuming all other countries in an `Other' category (thereby ensuring that no flow is lost). \highlight{Entries for which neither the importer nor the exporter have reported a value are assumed to be zero, and we constrain the cost matrix to be maximal on these entries. Entries for which only the exporter or the importer have reported a value (but not the other) are presumed missing in the respective table, and are masked in the loss function. With this approach, on average about 20\% of entries are masked in the transport plan (see fig.~\ref{fig:SI:masking} in the SI). Entries for which all reported values are missing populate the zero-flow edge set $\mathcal{S}$.}

\subsection{\highlight{Uncertainty quantification}}
Uncertainty on the estimates stems from two sources: one, the degree to which the minimizer of the inverse problem \eqref{eq:loss_function} is ill-defined (i.e. the number of possible cost matrices that all fit the problem equally well), and two, the uncertainty on the transport plans themselves. To address the first, we use an \emph{ensemble training approach} \cite{Gaskin_2023,Gaskin_2024} and train a family $\{u_k \}$ of 10 neural networks for each commodity in parallel. Even though the inverse problem is theoretically well-posed when the forward problem \eqref{eq:OT_entropy_regularised} admits a solution, in the case of the FAO data we are solving a minimisation problem. To incorporate the uncertainty on the transport plans, we pass random samples of $\bm{T}$ through each of the trained neural networks $u_k$. These samples are obtained by selecting either $T_{ij}^\text{E}$ or $T_{ij}^\text{I}$ uniformly at random for each entry of the transport plan, and passing this sample through each neural network. Repeating this $n$ times gives $n$ samples of $\bm{C}$, and inserting each estimate of $\bm{C}$ into Sinkhorn's algorithm gives $n$ estimated transport plans $\bm{\hat{T}}$. We generate a total of $n=100$ samples for each year and neural network (see also fig.~\ref{fig:SI:UQ} in the SI). The uncertainty estimates obtained by our method then provide an indication of how strongly a given set of trade flows informs the underlying cost.

\vspace{5mm}
\hrule
\vspace{1mm}
\subsection{Code and data availability}
All code and data is available at \url{https://github.com/ThGaskin/inverse-optimal-transport}. Instructions for running the model are given in the README.

\subsection{Acknowledgements}
TG, AD, and MTW acknowledge partial support by the EPSRC grant EP/X010503/1.  \highlight{AD also acknowledges support of a Royal Society APEX award (APX/R1/180133), under which this work was initiated.}

\vspace{5mm}
\hrule
\vspace{1mm}

\printbibliography

\appendix
\onecolumn
\renewcommand\thesection{S\arabic{section}}
\renewcommand\thetable{S\arabic{table}}
\setcounter{section}{0}
\setcounter{table}{0}
{\sffamily \noindent \huge \bfseries Supporting Information}
\vspace{2cm}
\renewcommand\thefigure{S\arabic{figure}}
\setcounter{figure}{0}
\phantomsection\addcontentsline{toc}{section}{Supporting Information}

\section*{\highlight{Validation of the approach on synthetic data}}
\begin{figure*}[b!]
\begin{adjustbox}{minipage=\textwidth-20pt}
\includegraphics[width=1\linewidth]{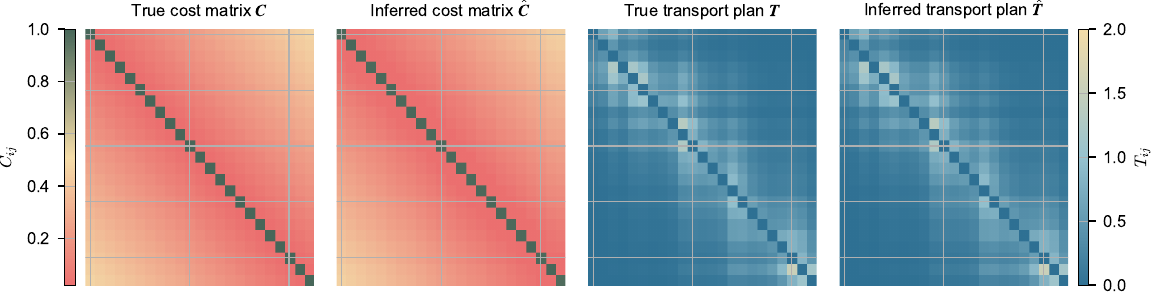}
\caption
{Inferring the cost matrix on noiseless, gapless, synthetic data. When requiring $C_{ij} \in [0, 1]$ and $u(0) = 1 = C_{max}$, the cost matrix $\bm{C}$ can be uniquely inferred from observations of the transport plan.}
\label{fig:SI:validation}
\end{adjustbox}
\end{figure*}

We first validate the neural inference approach on a synthetic, noiseless, gapless transport plan $\bm{T}$ with marginals $\bm{\mu}$ and $\bm{\nu}$, generated from the forward OT model with entropy regularisation $\varepsilon=0.15$, and infer the cost matrix. The diagonal of the transport plan is set to $0$, and the cost constrained to $[0, 1]$ by using a sigmoid activation function on the output layer of the neural network. By requiring the neural network to map the diagonal of $\bm{T}$ to the maximum permissible cost of $1$, the cost is uniquely inferred from $\bm{T}$ (fig.~\ref{fig:SI:validation}). Training is performed using the loss function described in the main manuscript, eq.~\eqref{eq:loss_function}.

\begin{figure}[hpt!]
\begin{adjustbox}{minipage=\textwidth-20pt}
    \includegraphics[width=1\linewidth]{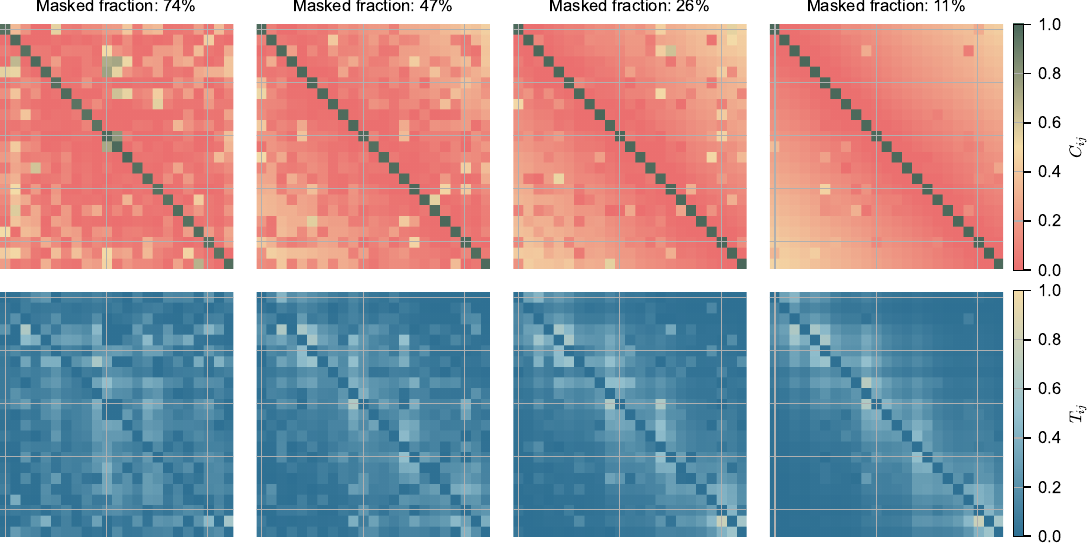}
    \vspace{0.1em}
    
    \includegraphics[width=1\linewidth]{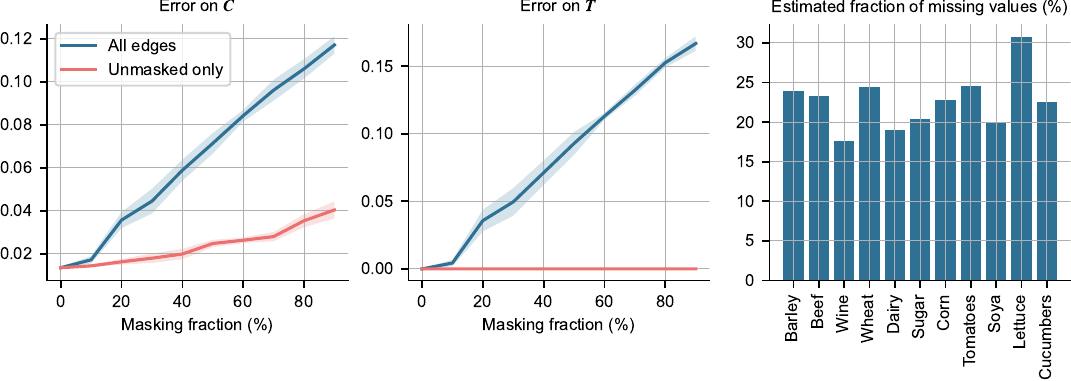}
    \caption[Performance on masked OT transport plans]{Top row: we infer the cost matrix on noiseless synthetic data with different fractions of the transport plan masked. Middle row: inferred transport plans. The ground truth is the same as in figure \ref{fig:SI:validation}. Bottom row: As the fraction of masked values increases, the average $L^1$ error on $\bm{T}$ and $\bm{C}$ increases (blue). However, on the non-masked values, errors remain significantly lower, and the prediction on the training values of $\bm{T}$ is independent of the proportion of masked values (red). Shown are mean and median values, as well as the standard deviation, over all entries of the respective matrices. Right: the estimated fraction of missing values in each FAO dataset. This is estimated by comparing the number of entries reported by one reporter (exporter/importer) but not the other.}
    \label{fig:SI:masking}
\end{adjustbox}
\end{figure}

\begin{figure}[ht!]
\begin{adjustbox}{minipage=\textwidth-20pt}
    \includegraphics[width=1\linewidth]{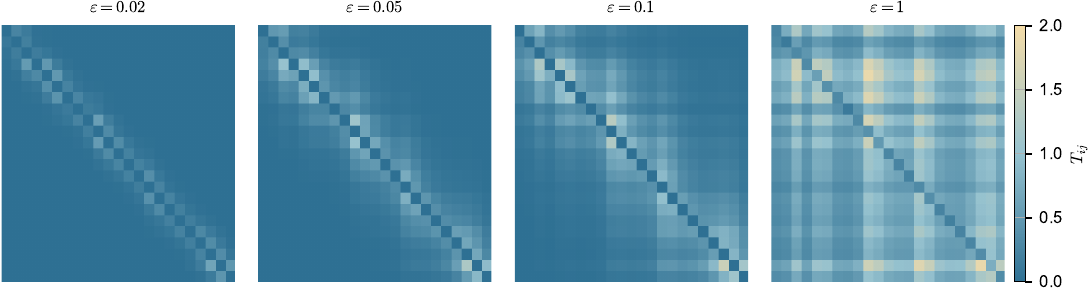}
    \vspace{0.1em}

    \includegraphics[width=1\linewidth]{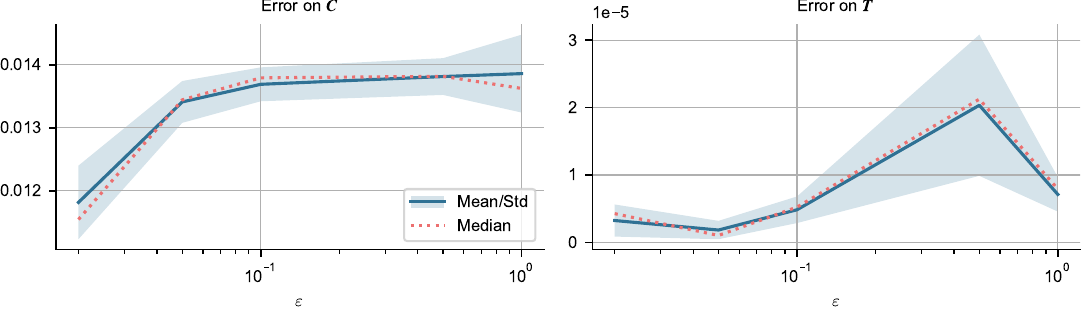}
    
    \caption[Effect of the OT entropy regulariser]{Effect of the entropy regulariser $\varepsilon$ on the inference. Top row: the transport plan for the cost matrix given in fig. \ref{fig:SI:validation} but diferent values of $\varepsilon$. Bottom row: accuracy on the inferred transport plan and cost matrix as a function of the regulariser $\varepsilon$. Each line is an average over 5 different seeds.}
    \label{fig:SI:epsilon}
\end{adjustbox}
\end{figure}

To analyse robustness with regard to missing values of $\bm{T}$, we mask a random number of entries in the transport plan, and re-infer the cost matrix $\bm{C}$, given the marginals of the unmasked transport plan. Results are shown in fig.~\ref{fig:SI:masking}. The method is robust for small amounts of missing data, since the number of gaps in each row and column will be small, thus constraining how the missing ``mass'' contained in the marginals can be distributed among the missing entries. The error on the non-masked values of $\bm{T}$ remains approximately constant (which is unsurprising), while the error on the entire cost matrix increases linearly with the masking fraction, though the error on the unmasked edges remains significantly smaller (red). For the FAO data, we estimate the number of missing data points from the number of entries that have exporter- but not importer-reported figures (or vice versa); on average, this gives a missing data fraction of around 20\% (see fig.~\ref{fig:SI:masking}).

Lastly, we analyse the effect of the choice of the entropy regulariser $\varepsilon$ on inference performance. $\varepsilon$ is a scaling parameter that determines how much small costs affect the transport plan (fig.~\ref{fig:SI:epsilon}), and can take any value in $[0, 1]$. One would like $\varepsilon$ to be as small as possible, since for $\varepsilon \to 0$ the inference procedure converges to classical OT. However, for small values of $\varepsilon$ Sinkhorn's algorithm becomes unstable, because as the entries of the initial guess $\exp(-\bm{C}/\varepsilon)$ go to $0$, the required scaling vectors need to grow exponentially to match the marginal constraints. This also causes the convergence rate of Sinkhorn's algorithm to slow significantly as $\varepsilon \to 0$, increasing the computational cost. We therefore choose a small value of $\varepsilon \approx 0.1$ that balances numerical stability and computational cost. The inference accuracy is independent of the choice of $\varepsilon$, see fig.~\ref{fig:SI:epsilon}.

\section*{\highlight{Uncertainty quantification}}
\begin{figure}[ht!]
\begin{adjustbox}{minipage=\textwidth-20pt}
    \includegraphics[width=1\linewidth]{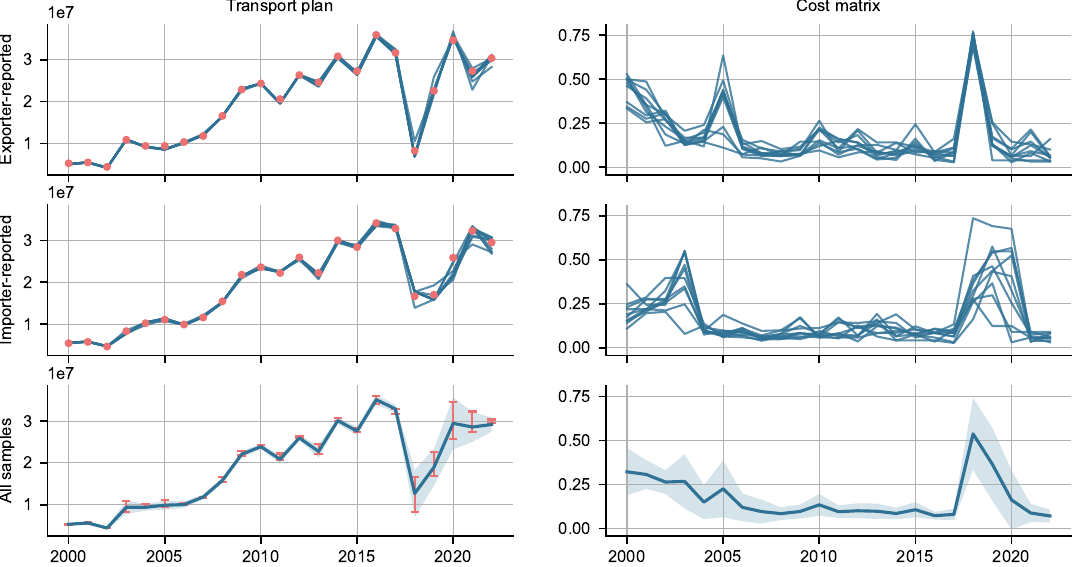}    
    \caption{Constructing the distribution on the cost matrix. A family of neural networks is trained on exporter-reported values (top) and importer-reported values (middle row), producing a distribution on the cost matrix due to the potential non-uniqueness of the minimizer. Bottom row: transport plan samples are generated by randomly mixing exporter- and importer-reported values and passing these through the neural network ensemble to produce the final distribution on $\bm{C}$ (right). }
    \label{fig:SI:UQ}

    \vspace{5mm}
    
    \includegraphics[width=1\linewidth]{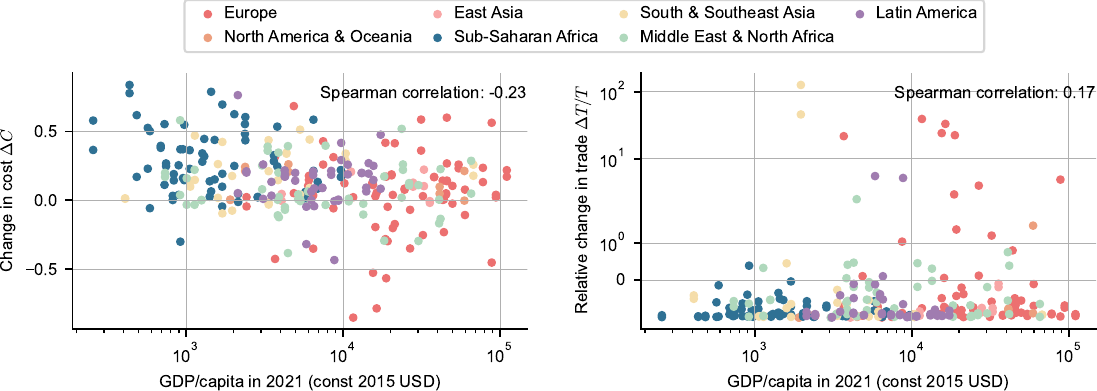}    
    \caption{Changes in trade costs (left) and relative change in trade volumes (right) with both Ukraine and Russia as a function of GDP per capita. Each dot represents the change for a single importing country, colour-coded by geographic region. Also shown are the Spearman correlation coefficients.}
    \label{fig:SI:Wheat_macrotrends_1}
\end{adjustbox}
\end{figure}

Fig.~\ref{fig:SI:UQ} illustrates our approach to uncertainty quantification. Shown is the US-China soya bean trade (fig.~\ref{fig:Soya_trade}{\cs{\textbf{C}}} in the main manuscript). The training data consists of exporter-reported values (top row) and importer-reported values (middle row). A family of neural networks produces a distribution of cost matrices that optimally reproduce the trade data (right column). The final distribution is generated by randomly sampling exporter- and importer-reported data points, passing them through the neural network ensemble, and averaging over the resulting cost matrices. The distribution over the cost matrices (bottom row, right) captures the uncertainty on the trade data (bottom row, left).

\begin{figure*}[ht!]
\begin{adjustbox}{minipage=\textwidth-20pt}
\begin{minipage}{\textwidth}
	\cs{\textbf{A} Russian wheat exports in metric tons, 2021 (left) and 2022}
\end{minipage}

\begin{minipage}{\textwidth}
	\includegraphics[width=\textwidth]{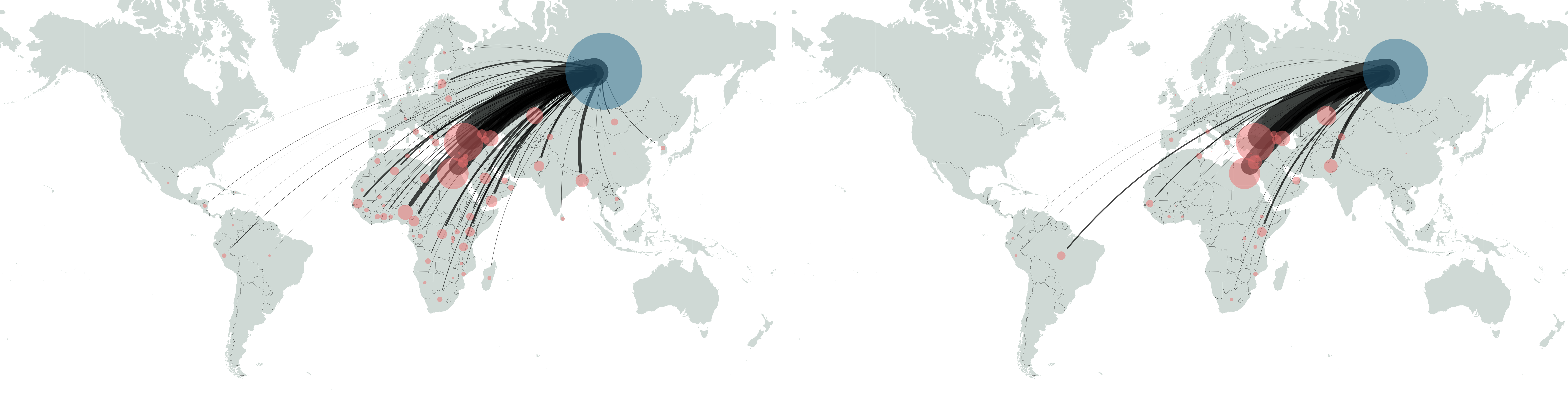}
\end{minipage}

\vspace{2mm}
\begin{minipage}[h]{\textwidth}
	\flushleft \cs{\textbf{B} Percent change in trade volume (left) and absolute change in cost}
\end{minipage}
\begin{minipage}{\textwidth}
	\includegraphics[width=\textwidth]{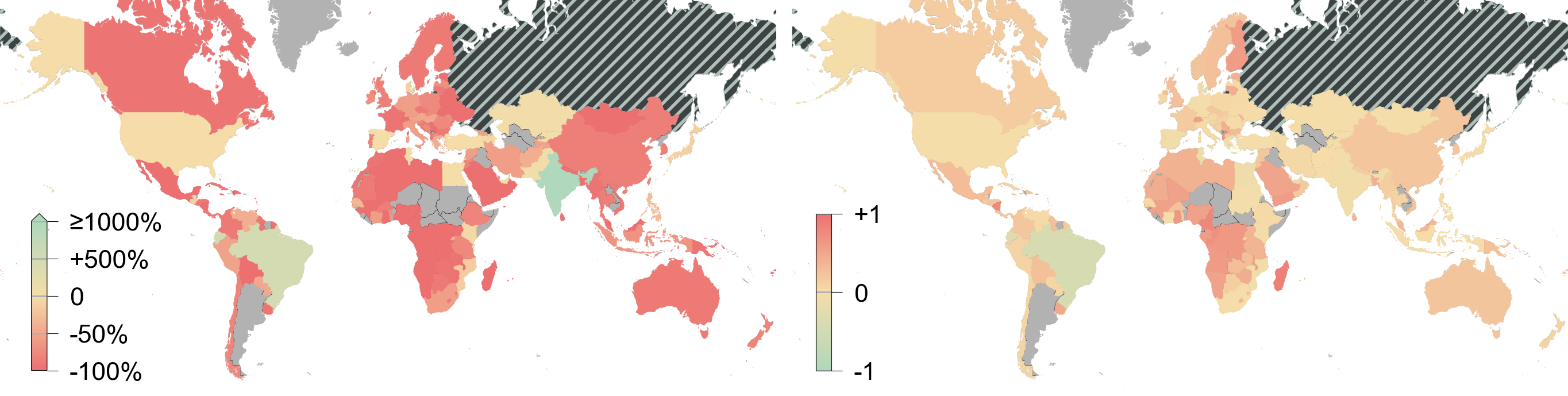}
\end{minipage}

\vspace{2mm}
\begin{minipage}{\textwidth}
	\cs{\textbf{C} Change in trade volume and cost, selected countries}
\end{minipage}
\begin{minipage}{\textwidth}
	\includegraphics[width=\textwidth]{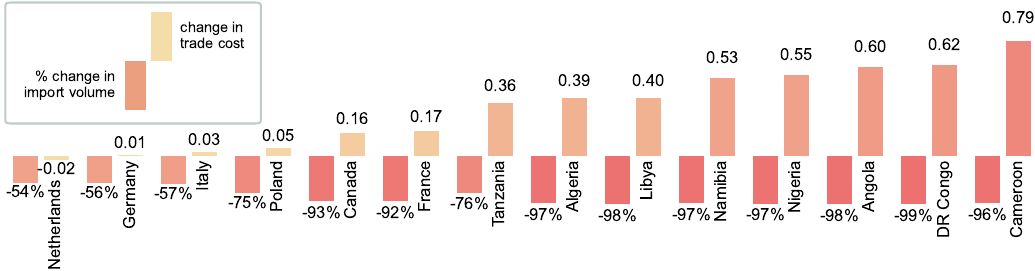}
\end{minipage}

\caption{The same plots as in figure \ref{fig:Ukraine_wheat} with Russia as the exporting partner.}
\label{fig:SI:Russia_wheat}
\end{adjustbox}
\end{figure*}

\section*{\highlight{Case study I: the impact of the Ukrainian war on wheat trade}}
\begin{figure}[ht!]
\begin{adjustbox}{minipage=\textwidth-20pt}
    \includegraphics[width=1\linewidth]{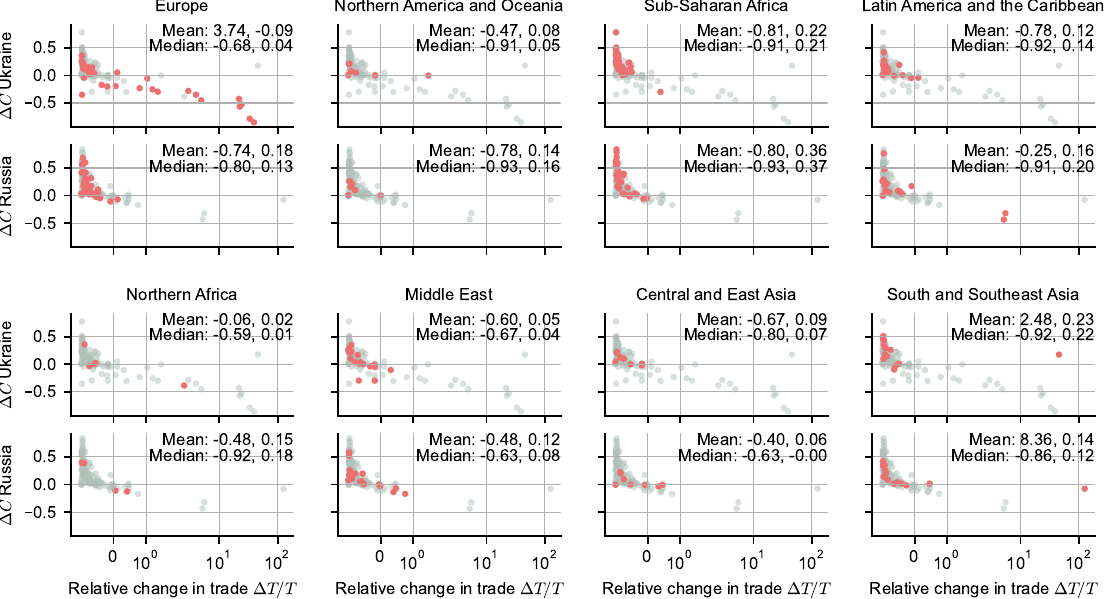}    
    \caption{Change in trade cost (y-axis) as a function of the relative change in trade volume (x-axis) for each region. The top rows show the change for Ukraine as the exporter, the bottom rows show the change for Russia as the exporter. Also given are the mean and median values for $\Delta T /T$ and $\Delta C$, respectively.}
    \label{fig:SI:Wheat_macrotrends_2}
    
    \vspace{1cm}
    
    \includegraphics[width=1\linewidth]{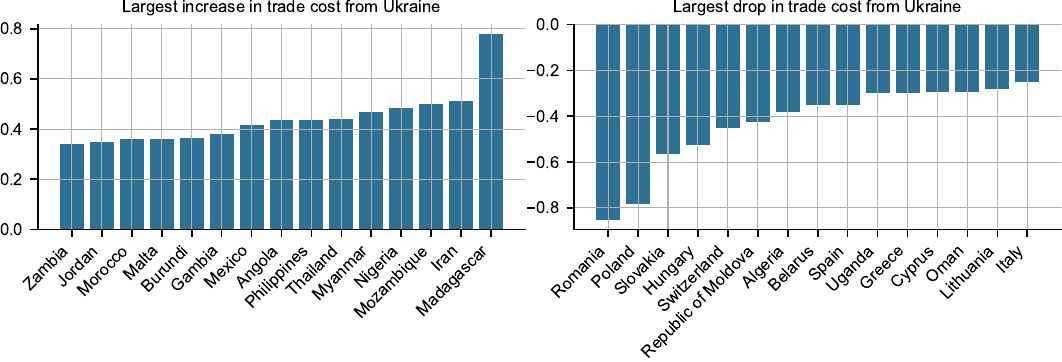}    
    \caption{The top 15 countries that experienced the largest increases (left) and largest drops (right) in trade costs with Ukraine, 2021--2022.}
    \label{fig:SI:Wheat_macrotrends_3}
\end{adjustbox}
\end{figure}

Here we show additional plots pertinent to our analysis of the global wheat market's response to the war in Ukraine. Figure \ref{fig:SI:Russia_wheat} is the equivalent of figure \ref{fig:Ukraine_wheat} in the main manuscript, but with Russia as the exporting partner. We see a similar trend of the increase in cost being disproportionately borne by the Global South. We substantiate this claim more thoroughly in figs. \ref{fig:SI:Wheat_macrotrends_1}, \ref{fig:SI:Wheat_macrotrends_2}, and \ref{fig:SI:Wheat_macrotrends_3}, which show changes in trade costs and trade volumes for each geographic region, as well as as a function of per-capita GDP. Figure \ref{fig:SI:Wheat_macrotrends_1} shows the change in trade costs (left) and relative change in trade volumes (right) from both Ukraine and Russia as a function of per-capita GDP, colour-coded by region. We see that European countries appear overrepresented among the countries that saw the largest drops in cost, while Sub-Saharan Africa is overrepresented among those with the largest increases in cost. The Spearman coefficient between change in cost and GDP/capita is $-0.23$, again indicating that wealthier countries suffered less from the impact of the war. This tallies with the fact that relative change in imports correlates positively with GDP/capita (right). 

Figure \ref{fig:SI:Wheat_macrotrends_2} shows the relative change in trade volumes $\Delta T/T$ from 2021--2022 and the associated drop in trade costs for each country, disaggregated by region. Europe, for instance, experienced a median 68\% decrease in trade volume and a 4\% increase in trade costs from Ukraine; North America saw a 91\% decrease in trade volume and a 5\% increase in trade costs. This contrasts with Sub-Saharan Africa (91\% decrease in trade, 21\% increase in costs), Latin America (92\% decrease in trade, 14\% increase in costs), or South and Southeast Asia (86\% decrease in trade, 22\% increase in costs). A similar, though slightly less pronounced pattern holds for Russian imports (bottom rows of each panel): here again, Sub-Saharan Africa saw the largest increase in costs, despite seeing a drop in trade comparable to that of Northern America and Oceania.

Figure \ref{fig:SI:Wheat_macrotrends_3} shows the 15 countries that experienced the largest drops and largest increases in trade costs with Ukraine between 2021--2022.

\vspace*{4in}

\begin{table}[t!]
\flushleft
\selectfont\tablefont
\begin{tblr}{
    width=\textwidth, 
    colspec={X[1.1,l]|X[c]X[c]X[c]X[c]X[c]X[c]X[c]X[c]X[c]}, 
    row{odd}={bg=box_color},row{1}={bg=box_color_dark},column{1}={bg=box_color_dark}
    }
 & $\lambda_1$ & $\lambda_2$ & $\lambda_3$ & $\lambda_4$ & $\lambda_5$ & $\lambda_6$ & $\lambda_7$ & $\lambda_8$ & $\lambda_9$ \\
\hline
Barley & 0.61 (0.071) & 0.085 (0.017) & -0.43 (0.017) & 1.2 (0.093) & -0.44 (0.12) & 0.72 (0.084) & 1.3 (0.1) & 0.66 (0.056) & -0.0058 (0.023) \\
Beef & 0.34 (0.041) & 0.18 (0.015) & -0.19 (0.017) & 1.7 (0.057) & -0.29 (0.068) & 0.78 (0.079) & 1.5 (0.068) & 0.55 (0.035) & 0.046 (0.017) \\
Corn & 0.86 (0.042) & 0.11 (0.012) & -0.49 (0.018) & 1.6 (0.077) & -0.67 (0.13) & 0.44 (0.081) & 1.4 (0.064) & 0.72 (0.059) & -0.0058 (0.023) \\
Cucumbers & 0.24 (0.047) & 0.27 (0.024) & -0.57 (0.023) & 1.8 (0.086) & -0.68 (0.1) & 1.6 (0.082) & 3.4 (0.21) & 0.46 (0.047) & 0.003 (0.012) \\
Dairy products$^\star$ & 0.58 (0.064) & 0.15 (0.0089) & -0.19 (0.0099) & 1.5 (0.055) & -0.11 (0.056) & 0.79 (0.042) & 1.9 (0.044) & 0.38 (0.026) & -0.066 (0.015) \\
Lettuce & 0.37 (0.044) & 0.13 (0.015) & -0.34 (0.014) & 1.9 (0.064) & -0.3 (0.085) & 0.7 (0.059) & 2.7 (0.12) & 0.31 (0.039) & -0.029 (0.014) \\
Soya & 0.81 (0.073) & 0.21 (0.019) & -0.069 (0.023) & 1.8 (0.13) & -0.84 (0.13) & 1.1 (0.16) & 1.2 (0.091) & 0.71 (0.039) & -0.053 (0.034) \\
Sugar products$^\dagger$ & 0.33 (0.033) & 0.094 (0.0077) & -0.33 (0.0086) & 2 (0.046) & 0.25 (0.051) & 0.46 (0.043) & 1.9 (0.038) & 0.58 (0.029) & 0.014 (0.012) \\
Tomatoes & 0.28 (0.063) & 0.15 (0.017) & -0.5 (0.021) & 1.8 (0.081) & -0.69 (0.11) & 2.2 (0.2) & 3 (0.14) & 0.31 (0.032) & 0.019 (0.018) \\
Wheat & 0.74 (0.042) & 0.091 (0.0073) & -0.41 (0.011) & 1.6 (0.055) & 0.59 (0.071) & 0.63 (0.048) & 1.2 (0.035) & 0.56 (0.029) & -0.012 (0.013) \\
Wine & 0.39 (0.077) & 0.18 (0.013) & -0.053 (0.012) & 1 (0.06) & 0.33 (0.06) & 1.1 (0.046) & 0.88 (0.048) & 0.34 (0.027) & -0.024 (0.019)
\end{tblr}
\caption{Estimated coefficients for each of the covariates used in the gravity model \eqref{eq:Gravity_model}, for each commodity. Standard errors are given in parentheses. The high-dimensional coefficients $\kappa_i$, $\omega_j$, and $\alpha_t$ are not shown. \smallskip \\ {\fontsize{7}{8}\selectfont $^*$Dairy products comprise: butter, skim milk of cows, cheese, other dairy products. $^\dagger$Sugar products comprise: sugar, refined sugar, syrups, fructose, sugar confectionery.}}
\label{tab:Gravity_parameters}
\end{table}
\section*{Comparison with Gravity model}
\begin{figure*}[htp!]
\begin{adjustbox}{minipage=\textwidth-20pt}
\begin{minipage}{0.5\textwidth}
	\flushleft \cs{\textbf{E} Dairy products} \vspace{2mm}
\end{minipage}
\begin{minipage}{0.5\textwidth}
	\flushleft \cs{\textbf{F} Sugar products} \vspace{2mm}
\end{minipage}
\begin{minipage}{0.5\textwidth}
	\includegraphics[width=\textwidth]{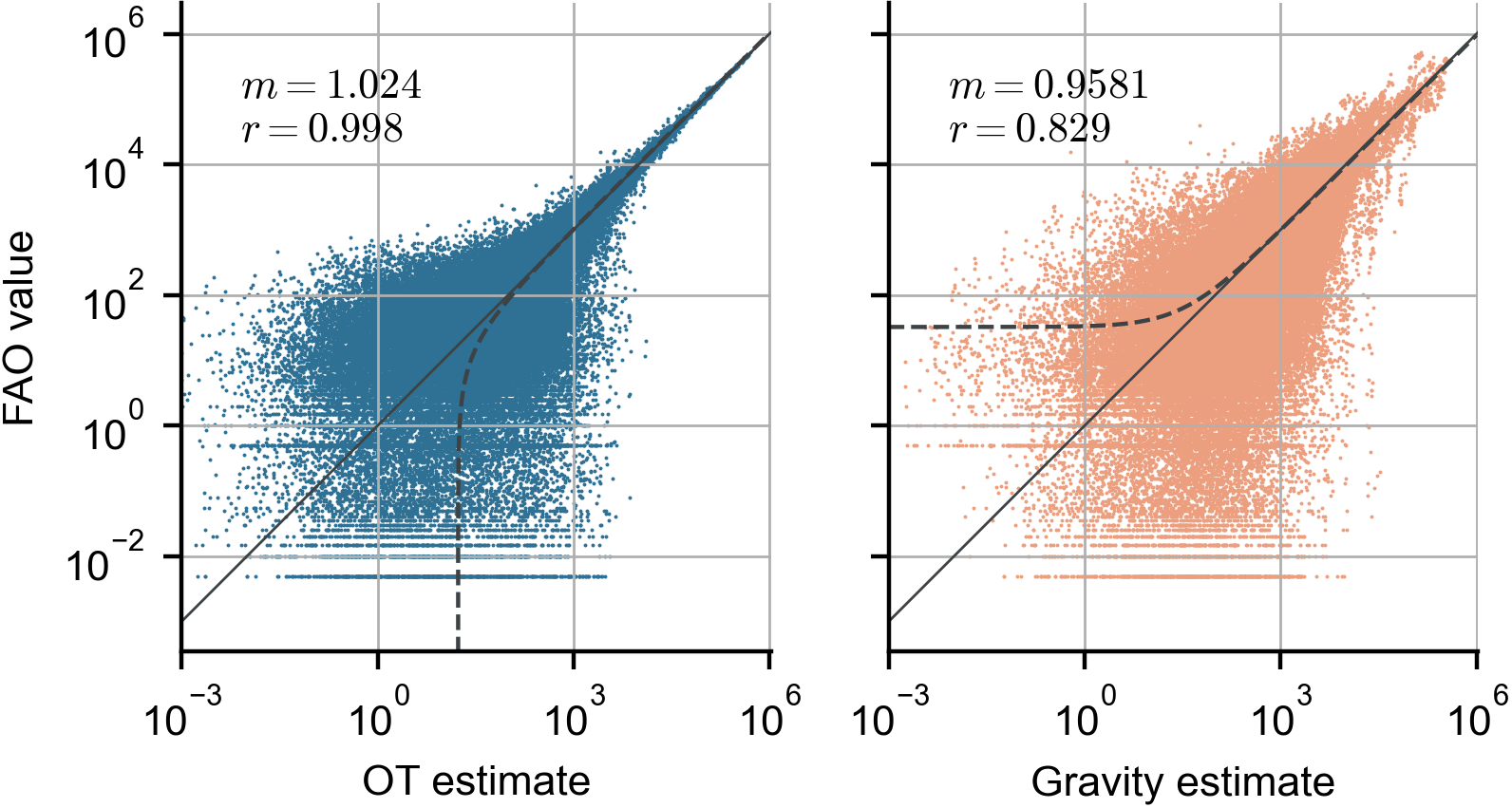}
\end{minipage}
\begin{minipage}{0.5\textwidth}
	\includegraphics[width=\textwidth]{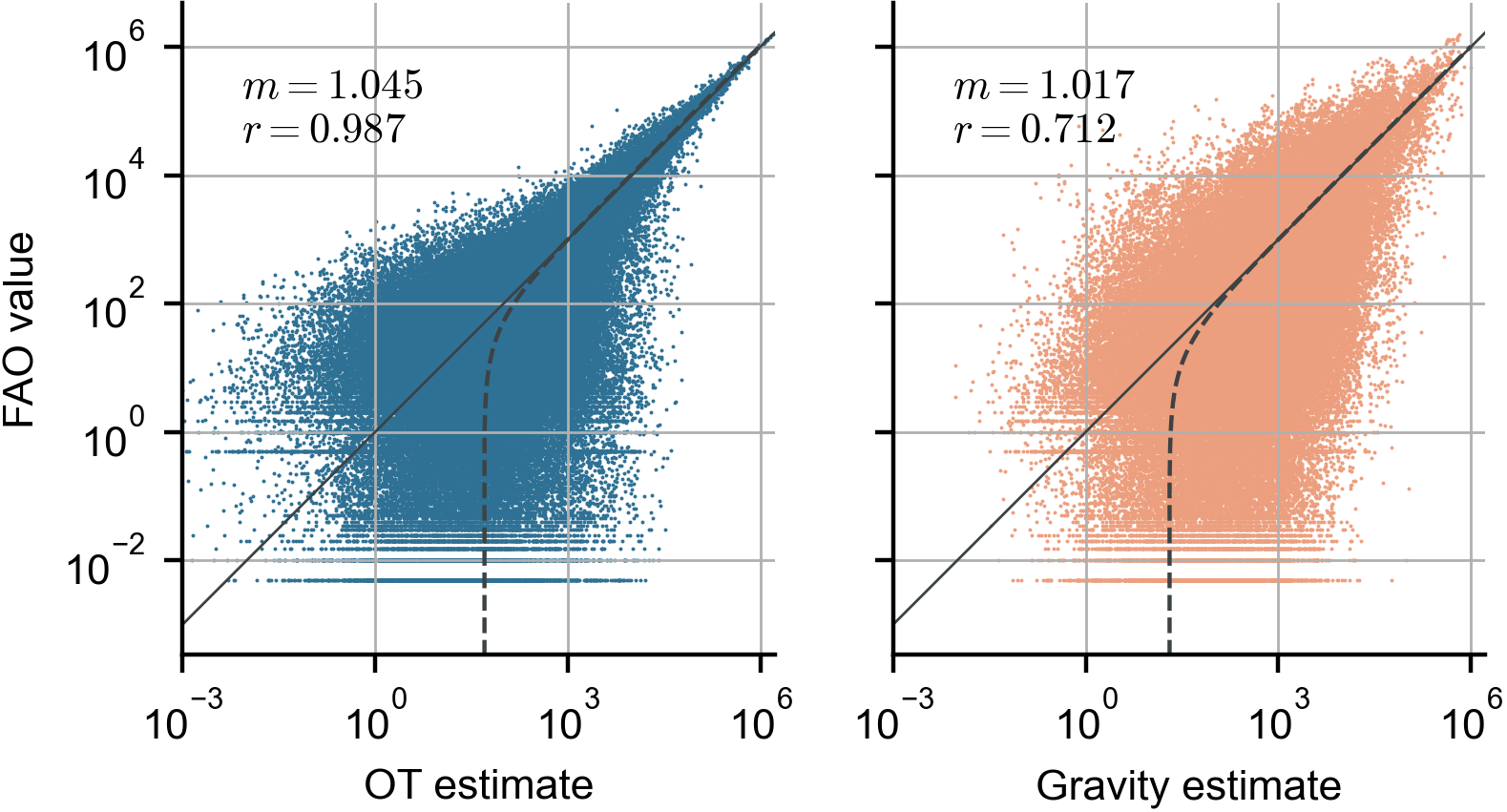}
\end{minipage}
\begin{minipage}{0.5\textwidth}
	\flushleft \cs{\textbf{G} Corn} \vspace{2mm}
\end{minipage}
\begin{minipage}{0.5\textwidth}
	\flushleft \cs{\textbf{H} Tomatoes} \vspace{2mm}
\end{minipage}
\begin{minipage}{0.5\textwidth}
	\includegraphics[width=\textwidth]{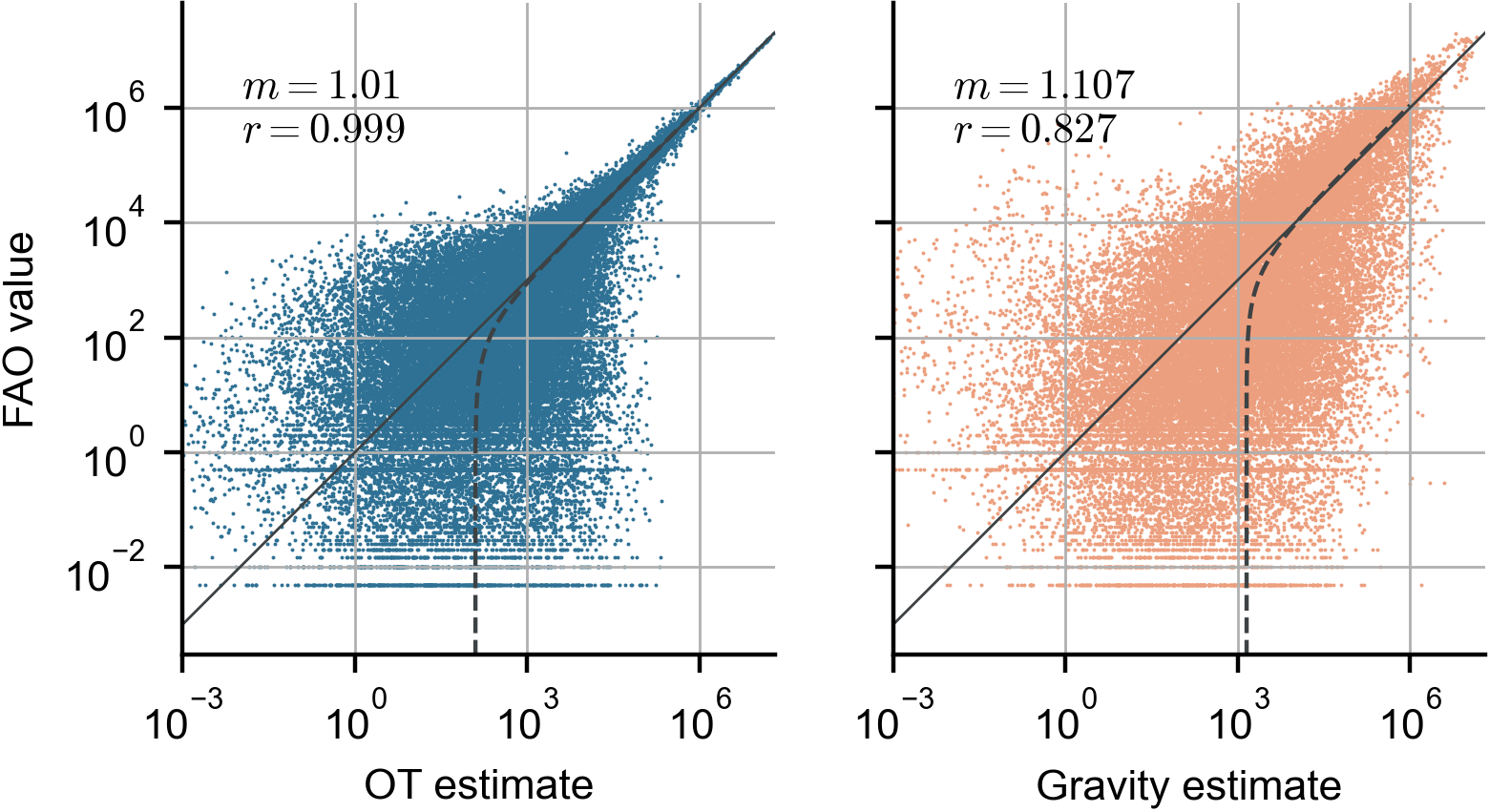}
\end{minipage}
\begin{minipage}{0.5\textwidth}
	\includegraphics[width=\textwidth]{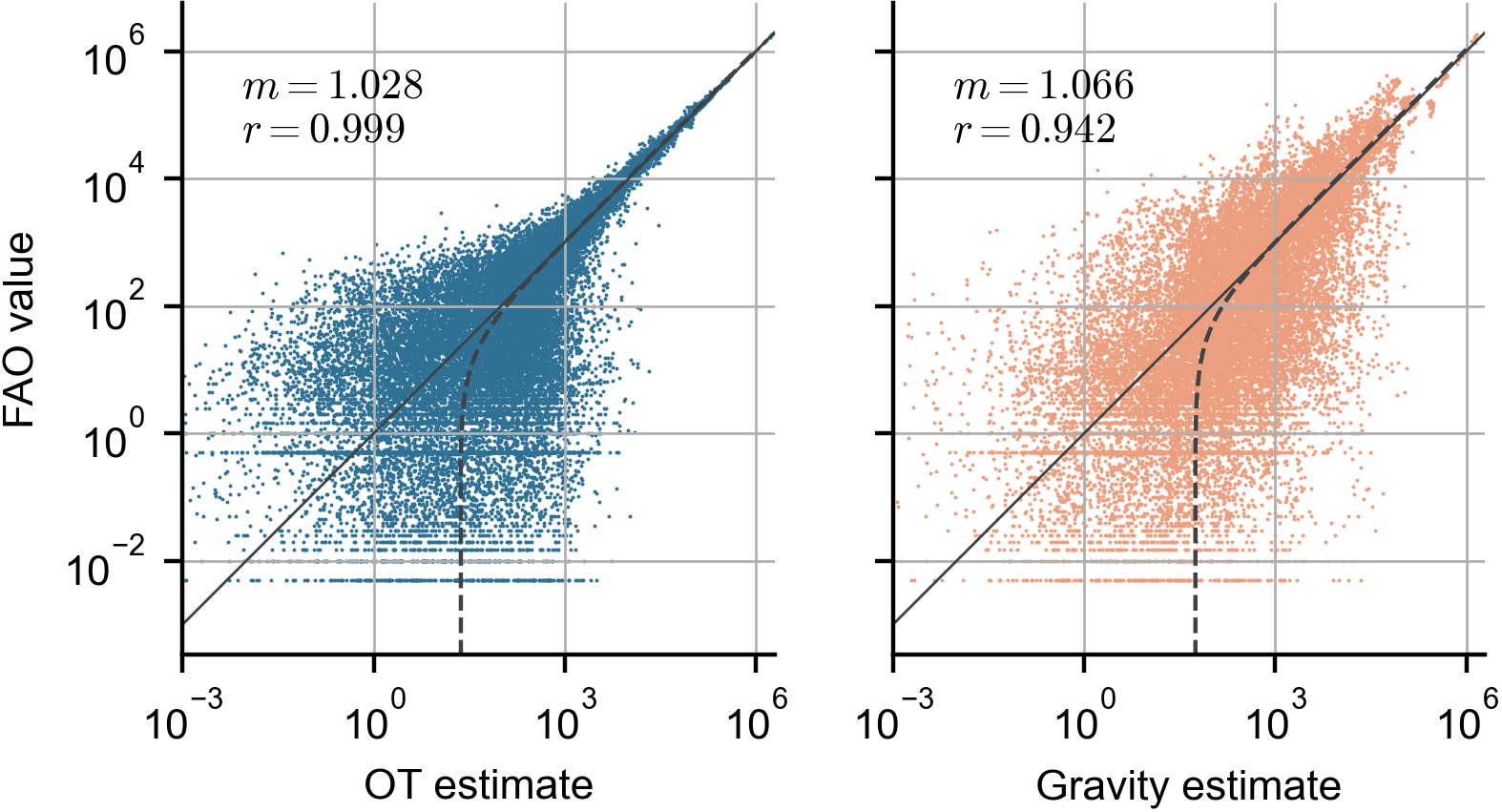}
\end{minipage}
\begin{minipage}{0.5\textwidth}
	\flushleft \cs{\textbf{I} Soya beans} \vspace{2mm}
\end{minipage}
\begin{minipage}{0.5\textwidth}
	\flushleft \cs{\textbf{J} Vegetables} \vspace{2mm}
\end{minipage}
\begin{minipage}{0.5\textwidth}
	\includegraphics[width=\textwidth]{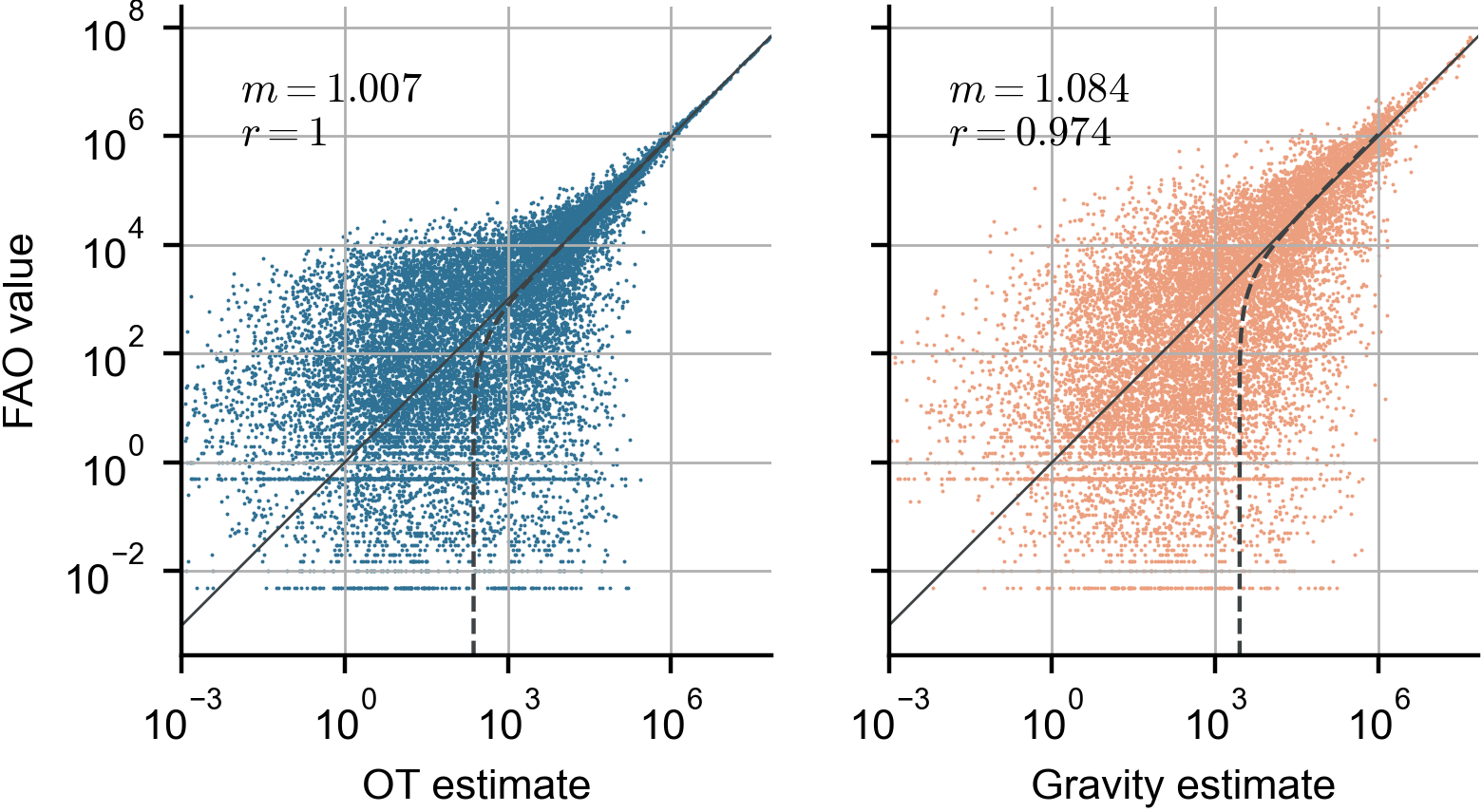}
\end{minipage}
\begin{minipage}{0.5\textwidth}
	\includegraphics[width=\textwidth]{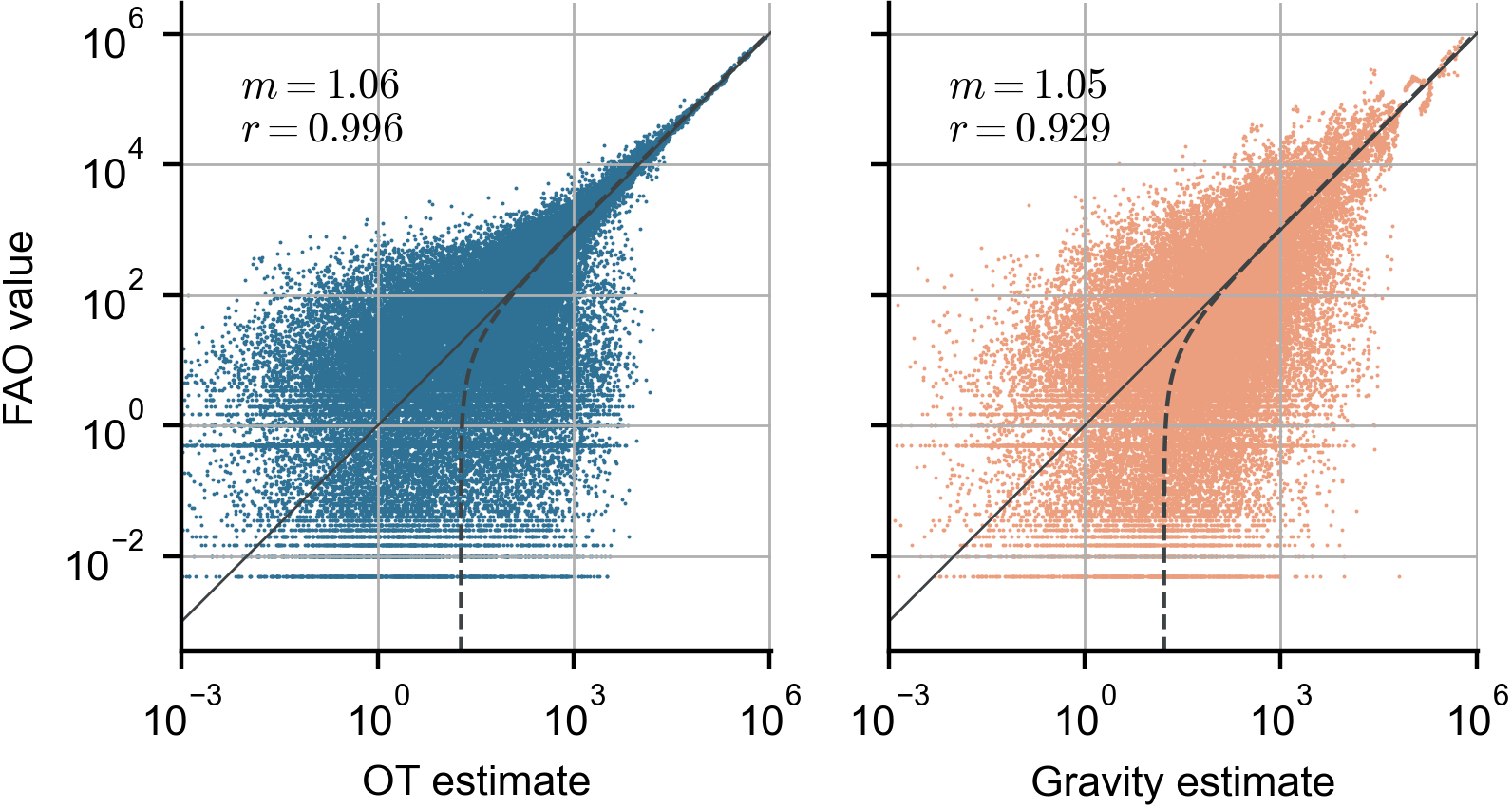}
\end{minipage}
\begin{minipage}{\textwidth}
	\flushleft \cs{\textbf{K} Cucumbers and gherkins} \vspace{2mm}
\end{minipage}
\begin{minipage}{0.5\textwidth}
	\includegraphics[width=\textwidth]{Figures/Regression_fits/Cucumbers.png}
\end{minipage}
\hspace{0.05\textwidth}
\begin{minipage}{0.45\textwidth}
    \caption{Comparison of the OT estimates (left, darkblue) and the Gravity estimates (eq.~\eqref{eq:Gravity_model}, for each commodity. The $y$-axis shows the true FAO value, while the $x$-axis shows the estimated value. The solid line is the diagonal $y=x$. Also shown are a linear fit (dashed line) as well the fitted slope $m$ and Pearson correlation $r$ of the fit.}
    \label{fig:SI:Gravity_comparisons_all}
\end{minipage}
\end{adjustbox}
\end{figure*}

\highlight{We consider the following gravity model specification in the spirit of traditional gravity estimations, where the covariates are based on \cite{Yotov_2017}:
\begin{align}
    T_{i,j,t,l} = & \exp (\kappa_{i,l} + \omega_{j,l} + \alpha_{t,l} + \lambda_1 \log O_{i,t,l} + \lambda_2 \log E_{j,t,l} + \lambda_3 \log d_{i,j} + \lambda_4 \mathrm{CNTG}_{i,j} \nonumber \\ & + \lambda_5 \mathrm{CNLY}_{i,j} + \lambda_6 \mathrm{LANG}_{i,j} + \lambda_7 \mathrm{RTA}_{i,j} + \lambda_8 \log \chi_{j,t,l} + \lambda_9 \log \mathrm{TRFF}_{i,j,t,l})\varepsilon_{i,j,t,l}.
    \label{eq:Gravity_model}
\end{align}
The covariates 5--9 are taken from the CEPII database \cite{CEPII}:
\begin{enumerate}
    \item $\kappa_{i,l}$ are the exporter-fixed effects,
    \item $\omega_{i,l}$ are the importer-fixed effects,
    \item $\alpha_{t,l}$ are the year-fixed effects,    
    \item $O_{i,t,l}$ is the total production output, in tonnes, of the exporter $i$ of product $l$ at time $t$ as given by the FAO,
    \item $E_{j,t,l}$ is the total consumption of the importer $j$ of product $l$ at time $t$, in tonnes,
    \item $d_{i,j}$ is the geodesic distance in km between the population centres (harmonic average) of countries $i$ and $j$ (\texttt{distw\_harmonic}),
    \item $\mathrm{CNTG}_{i,j}$ is the binary variable that indicates whether countries $i$ and $j$ share a land border (\texttt{contig}),
    \item $\mathrm{CNLY}_{i,j}$ is a binary variable indicating whether there ever existed colonial ties before 1948 between the two trading partners (\texttt{col\_dep\_ever}),
    \item $\mathrm{LANG}_{i,j}$ indicates whether countries $i$ and $j$ share an official or primary language (\texttt{comlang\_off}),
    \item $\mathrm{RTA}_{i,j,t}$ is a binary variable indicating whether there exists a regional trade agreement between countries $i$ and $j$ at time $t$ (\texttt{rta\_coverage}),
    \item $\chi_{j,t,l} = \sum_i d_{i,j} O_{i,t,l} / \sum_k O_{k,t,l}$ is the remoteness index of the importer,
    \item $\mathrm{TRFF}_{j,t,l}$ is the tariff applied by the importer $j$ to product $l$ at time $t$ in the absence of a trade agreement. We use the most favoured nation tariff (maximum duty) as given by the WTO \cite{WTO_Stats}: \texttt{MFN - Maximum duty by product groups}. 
\end{enumerate}
The remoteness index and the exporter/importer fixed effects account for the multilateral resistance terms \cite{Yotov_2017}. This gives a $L + m + n + 9$-dimensional regression problem for each commodity, where $L = 22, m, n$ denote the number of years, exporter countries, and importer countries in the dataset (note that the regressors only span the period until 2021).

\begin{figure*}[t!]
\begin{adjustbox}{minipage=\textwidth-20pt}
\begin{minipage}{0.46\textwidth}
	\flushleft \cs{\textbf{A} Estimation accuracy (RMSE)}
\end{minipage}
\begin{minipage}{0.54\textwidth}
	\flushleft \cs{\textbf{B} Estimation accuracy (RMSE in standard deviations)}
\end{minipage} \hfill \vspace{-2mm}

\begin{minipage}{\textwidth}
	\includegraphics[width=\textwidth]{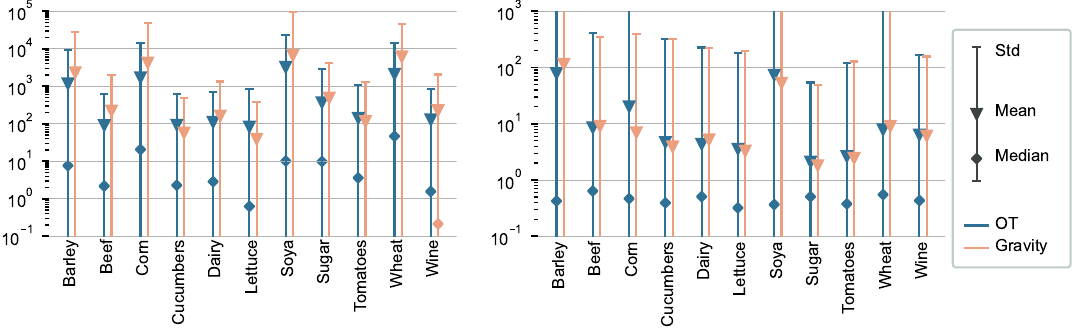}
\end{minipage}
\vspace{2mm}

\begin{minipage}{\textwidth}
    \cs{\textbf{C} Pearson Correlation coefficients}
\end{minipage}
\begin{minipage}{\textwidth}
	\includegraphics[width=\textwidth]{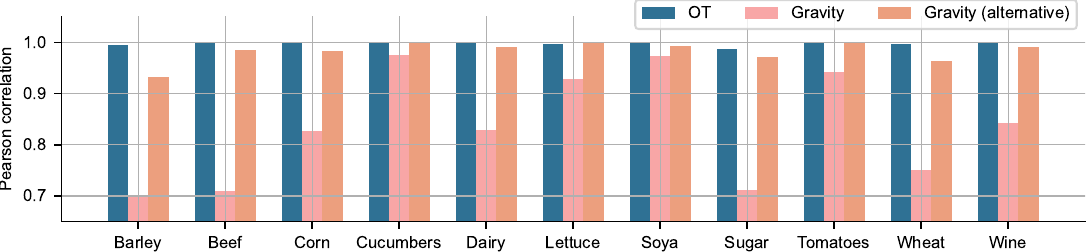}
\end{minipage}
\caption{{\cs{\textbf{A}--\textbf{B}}}~The same plots as in fig.~\ref{fig:Gravity_comp}{\cs{\textbf{C}--\textbf{D}}} using the alternative three-way gravity model \eqref{eq:Gravity_model_alt}. The $y$-axis scaling is preserved. Based on the median RMSE, the three-way gravity model outperforms OT, whereas OT performs better in terms of the mean. This reflects the greater weighting---and thus lower error---for large flows, as shown in Figure~\ref{fig:SI:Gravity_comparisons_all}. {\cs{\textbf{C}}} Pearson correlation coefficients between observed and estimated transport flows for OT and both gravity models (eqs.~\eqref{eq:Gravity_model} and \eqref{eq:Gravity_model_alt}).}
\label{fig:SI:Gravity_comp_alt}
\end{adjustbox}
\end{figure*}
We estimate the parameters of the gravity equations using the Poisson Pseudo Maximum Likelihood (PPML) estimator \cite{silva2006log,Yotov_2017}. Data of intra-national flows and zero trade flows are excluded for consistency with the optimal transport framework. We use the PPML with high-dimensional fixed effects implementation developed in \cite{Correia2020}. Table \ref{tab:Gravity_parameters} gives the estimated parameters for each commodity. Figure \ref{fig:SI:Gravity_comparisons_all} plots the estimated values $\hat{T}_{ij}$ against the reporter-averaged FAOStat values for both the OT and the gravity models. Also shown are a linear fit with slopes and Pearson coefficients indicated.

As an alternative specification, we use the three-way gravity model with exporter-time, importer-time, and pair (exporter-importer) fixed effects, absorbing exporter or importer level covariates, except for bilateral time variant trade policy variables, into the fixed effects, leading to:
\begin{equation}
    T_{i,j,t,l} = \exp(\gamma_{i,t,l} + \sigma_{j,t,l} + \beta_{i,j,l} + \lambda \mathrm{RTA}_{i,j,t})\varepsilon_{i,j,t,l}.
    \label{eq:Gravity_model_alt}
\end{equation}
This specification has the highest level of flexibility, using regressors only for time-dependent bilateral trade terms (participation in regional trade agreements). The country-time fixed effects capture both the time-dependent pull effect due to the size and the multilateral resistance terms, while the pair effect captures time-invariant bilateral effects \cite{Breinlich2021}. This alternative specification has $L \times (m + n) + m \times n + 1$ parameters (see fig.~\ref{fig:SI:Gravity_comp_alt}).}

\begin{figure*}[h!]
\begin{adjustbox}{minipage=\textwidth-20pt}
\begin{minipage}{0.5\textwidth}
	\flushleft \cs{\textbf{A} Barley}
\end{minipage}
\begin{minipage}{0.5\textwidth}
	\flushleft \cs{\textbf{B} Beef}
\end{minipage}
\hfill\vspace{-3mm}

\begin{minipage}{0.5\textwidth}
	\includegraphics[width=\textwidth]{Figures/Regression_fits/Barley.png}
\end{minipage}
\begin{minipage}{0.5\textwidth}
	\includegraphics[width=\textwidth]{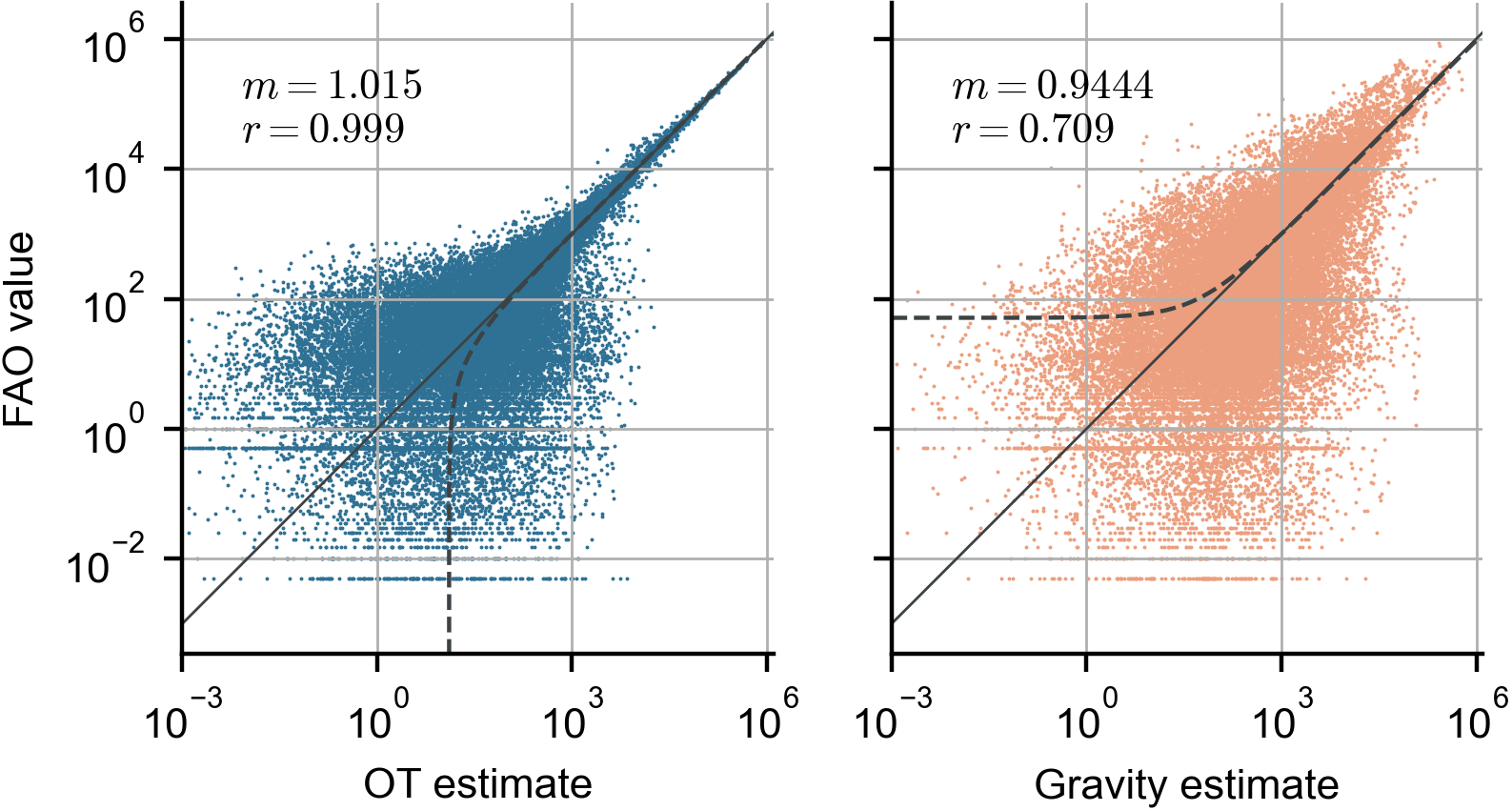}
\end{minipage}
\vspace{2mm}

\begin{minipage}{0.5\textwidth}
	\flushleft \cs{\textbf{C} Wine}
\end{minipage}
\begin{minipage}{0.5\textwidth}
	\flushleft \cs{\textbf{D} Wheat}
\end{minipage}
\hfill\vspace{-3mm}

\begin{minipage}{0.5\textwidth}
	\includegraphics[width=\textwidth]{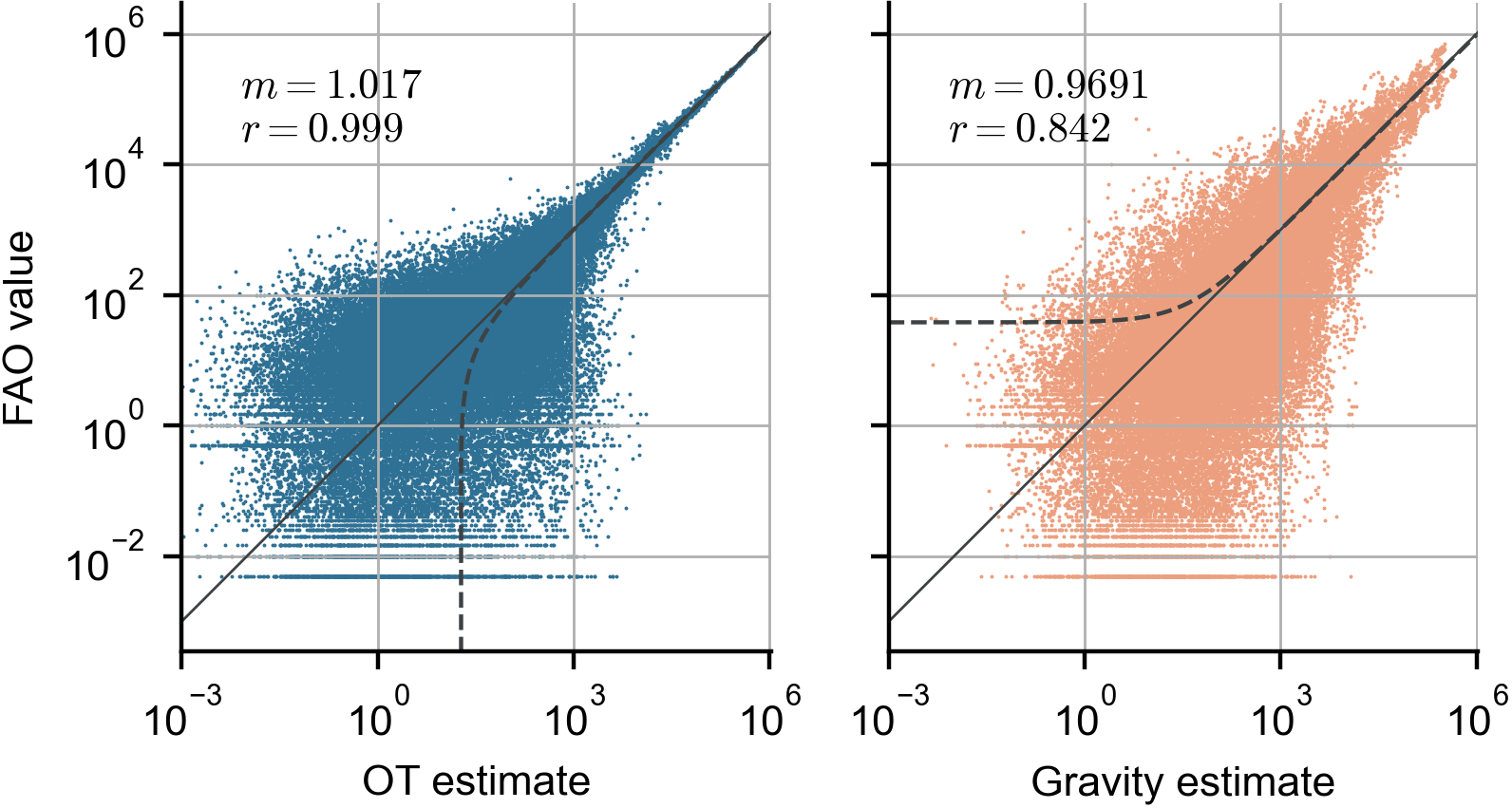}
\end{minipage}
\begin{minipage}{0.5\textwidth}
	\includegraphics[width=\textwidth]{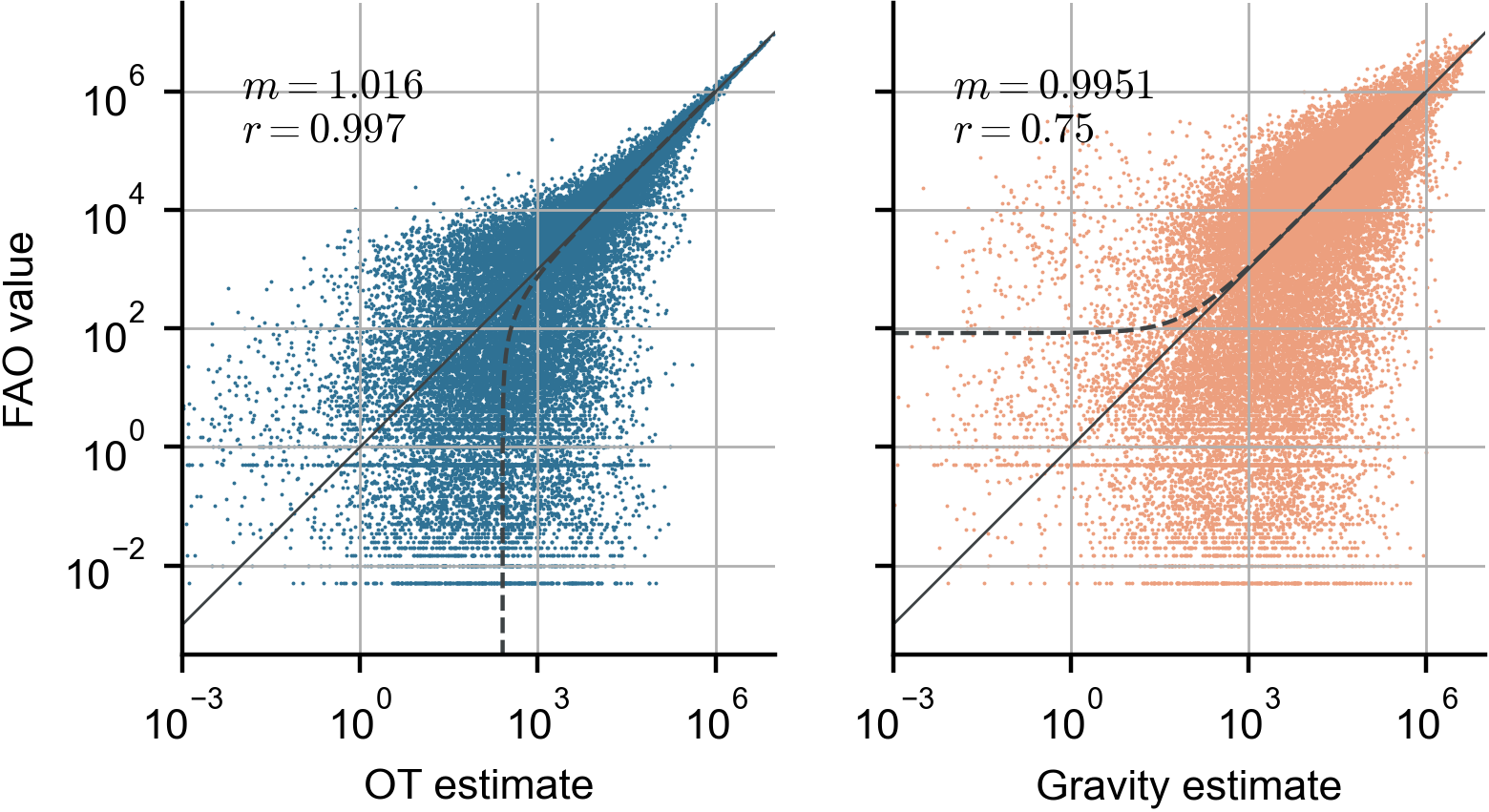}
\end{minipage}

\caption*{(continued on next page)}
\end{adjustbox}
\end{figure*}

\end{document}